\documentclass{jtemplates/CUP-JNL-DCE-CLEAN}
\pdfoutput=1

\usepackage{latexsym}
\usepackage{graphicx}
\usepackage{multicol,multirow}
\usepackage{rotating}
\usepackage{appendix}
\usepackage[authoryear]{natbib}
\usepackage{ifpdf}
\usepackage[T1]{fontenc}
\usepackage{times}
\usepackage{sourcesanspro}
\usepackage{textcomp}%
\usepackage{xcolor}%
\usepackage[draft]{hyperref}

\articletype{RESEARCH ARTICLE}
\jname{Data-Centric Engineering}
\jyear{2022}

\usepackage{moreverb,url}

\usepackage{nameref}

\usepackage{tikz}
\usetikzlibrary{arrows,calc,positioning,shapes}

\usepackage{csml}

\usepackage{hhline}
\usepackage{mathtools}
\usepackage{breqn}

\usepackage{floatrow}
\newfloatcommand{capbtabbox}{table}[][\FBwidth]
\usepackage{algorithm}
\usepackage[noend]{algpseudocode}

\newcommand\BibTeX{{\rmfamily B\kern-.05em \textsc{i\kern-.025em b}\kern-.08em
T\kern-.1667em\lower.7ex\hbox{E}\kern-.125emX}}

\begin{document}
	
	\begin{Frontmatter}
	
	\authormark{Febrianto, Butler, Girolami and Cirak}

	\title{Digital twinning of self-sensing structures using the statistical finite element method}

	\author[1,3]{Eky Febrianto}
	\author[2,3]{Liam Butler}  
	\author[1,3]{Mark Girolami}  
	\author*[1,3]{Fehmi Cirak} \email{f.cirak@eng.cam.ac.uk}
	
	\address[1]{\orgdiv{Department of Engineering}, \orgname{University of Cambridge}, \orgaddress{\street{Trumpington Street}, \state{Cambridge}, \postcode{CB2~1PZ}, \country{UK}}}
	\address[2]{\orgdiv{Lassonde School of Engineering}, \orgname{York University}, \orgaddress{ \state{Toronto}, \postcode{ON M3J 1P3}, \country{Canada}}}
	\address[3]{\orgname{The Alan Turing Institute}, \orgaddress{\street{96 Euston Road}, \state{London}, \postcode{NW1 2DB}, \country{UK}}}



\abstract{
The monitoring of  infrastructure assets using sensor networks is becoming increasingly prevalent. A digital twin in the form of a finite element model,  as commonly used in design and construction, can help make sense of the copious amount of collected sensor data.  This paper demonstrates the application of the statistical finite element method (statFEM), which provides a principled means of synthesising data and physics-based models, in developing a digital twin of a self-sensing structure. As a case study, an instrumented steel railway bridge of $27.34~\text{m}$ length located along the West Coast Mainline near Staffordshire in the UK is considered. Using strain data captured from fibre Bragg grating (FBG) sensors at 108 locations along the bridge superstructure, statFEM can predict the `true' system response while taking into account the uncertainties in sensor readings, applied loading and finite element model misspecification errors. Longitudinal strain distributions along the two main I-beams are both measured and modelled during the passage of a passenger train. The statFEM digital twin is able to generate reasonable strain distribution predictions at locations where no measurement data is available, including at several points along the main I-beams and on structural elements on which sensors are not even installed.  The implications for long-term structural health monitoring and assessment include optimisation of sensor placement and performing more reliable what-if analyses at locations and under loading scenarios for which no  measurement data is available.
}

	\keywords{digital twin, structural health monitoring, statistical finite element method, Bayesian learning, physics-informed machine learning, FBG sensing}
	
	\begin{policy}[Impact Statement]
	Engineering structures, such as bridges, tunnels, power plants or buildings, are vital to the functioning of society. So far, they have been designed, built and maintained primarily with the aid of computational finite element models that rely on numerous empirical assumptions and codified safety factors. The predictions of such models often bear little resemblance to the behaviour of the actual structure. Lately, advances in infrastructure monitoring using sensor networks is providing an unprecedented amount of data from structures in operation. The recently proposed statistical finite element method boosts the predictive accuracy of models by synthesising their output with sensor data. In this work, the statistical finite element method is used in developing a digital twin of an operational railway bridge. The obtained digital twin enables the autonomous continuous monitoring and condition assessment of the bridge.
	\end{policy}

	\end{Frontmatter}


	%
\section{Introduction \label{sec:introduction}}
%
%

Designers, engineers and maintenance managers commonly rely on finite element~(FE) models to help better understand the current and future behaviour of infrastructure assets. FE models are used both to facilitate the design of new structures and to help assess structures that have already been in operation for many years. The level of sophistication of the  FE models depends on their intended purpose, the structure modelled, and the resources available (i.e., technical skill, time, computational power, etc.). As with any other model, FE models rely on assumptions and parameters, which are inevitably subject to uncertainties, both epistemic and aleatoric,  and require sufficient degree of experience to verify and interpret the predicted results~\citep{oden2010computerI,lau2018role}.

More recently, a trend toward monitoring of infrastructure assets using sensor networks has led to the development of new data-driven approaches to assess the performance of structures; see the reviews~\citep{abdulkarem2020wireless,lynch2007overview,brownjohn2007structural}. Sensing data can, for instance, provide performance information to help decide when to inspect, repair or decommission a structure and to accelerate or improve prototyping of new construction techniques and materials~\citep{frangopol2016life}. There has been a wide array of field studies that have implemented sensing in bridges, tunnel linings, high-rise buildings and dams.  For instance,~\cite{deBattista2017monitoring} instrumented several columns and walls on every floor in a 50-storey residential tower in central London.  Using distributed fibre optic sensors (DFOS), they measured the axial shortening of the various columns and walls as the building was being erected.  Recently,~\cite{di2021fatigue} reported a study in which electrical-based strain gauges were installed on a $450~\text{metre}$ steel tied-arch bridge with an orthotropic steel deck.  These types of bridge structures are particularly susceptible to fatigue damage and, using the data captured via the sensors, they were able to identify several fatigue-vulnerable locations along the structure.  Another unique application of structural health monitoring of a $100~\text{metre}$ segment of tunnel lining at the CERN (European Council for Nuclear Research) was conducted by~\cite{diMurro2016distributed}.  Localised tension and compression cracks were observed in the tunnel and prompted instrumentation and monitoring using DFOS.  In light of data captured over ten months they reported strain levels in the tunnel, which ultimately proved to be insignificant.

With the increase in structural health monitoring campaigns, the infrastructure engineering sector has been inundated with large quantities of data whose value is still being realised. In particular, this emerging sensing paradigm has resulted in the development of new data-driven strategies to model infrastructure performance; see the insightful  reviews~\cite{huang2019state,wu2020data}.  Unlike FE modelling, data-driven approaches utilising monitoring data provide information about the actual operational performance of a particular infrastructure asset. The modelling effort associated with data-driven techniques is lower than in FE modelling. However, even the most advanced data-driven models, which are often based on artificial neural networks, Gaussian process regression or Bayesian model updating, require a copious amount of \emph{real-world} training data to make predictions. This challenge has prompted research in the areas of model updating and system identification which have attempted to use both FE and data-driven models to provide better predictions of structural behaviour, see e.g. ~\cite{malekzadeh2015performance, pasquier2016iterative,tsialiamanis2021generative}.  However, the current modelling approaches seem to be unable to synthesise measurement data with uncertainties and predictions from inherently misspecified FE models in a manner that allows for the generation of continuous predictions as new measurement data becomes available.  This ability is, however, critical to the realisation of digital twins~\citep{rasheed2020digital,worden2020digital}. In its broadest sense, as defined in a UK governmental report, a digital twin is a realistic digital representation of assets, processes or systems in the built or natural environment~\citep{bolton2018gemini}.  

 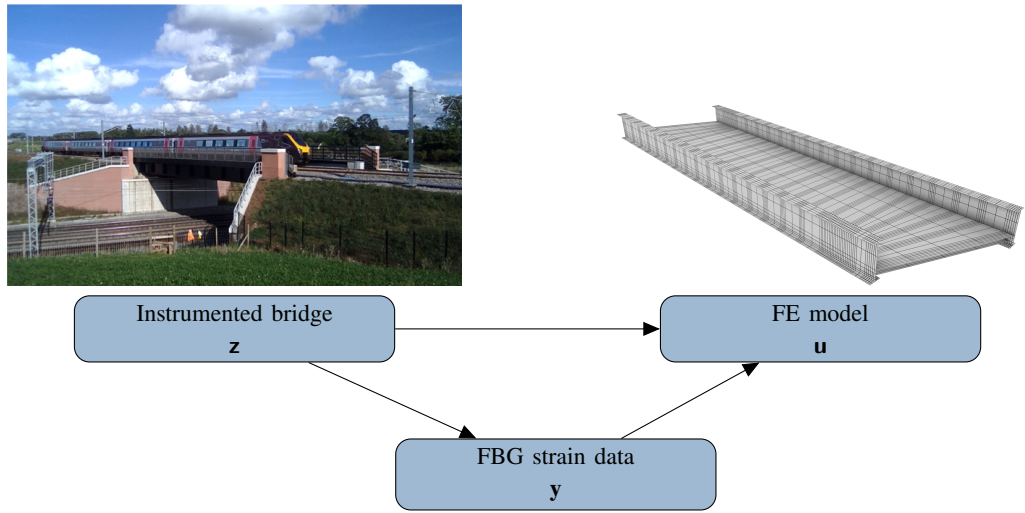
\begin{figure*}[h!]
	\centering
	\begin{tikzpicture}[
	->, >=triangle 45, shorten >=0pt, auto,node distance=65mm, 
	main node/.style={rectangle, rounded corners =2mm, fill=csmlLightBlue, draw=black, minimum width = 30mm, minimum height = 7.5mm, align=center, text width=40mm},
	pic node/.style={fill=none, draw=none, minimum width = 30mm, minimum height = 7.5mm, align=center, text width=60mm}, 
	scale=.5] 
	\node[main node]    (1)             {\small Instrumented bridge \\ $\ary z$};
	 
	\node[main node]    (2)            [right= 35 mm of 1]     {\small FE model \\ $\ary u$}; 
	\node[main node]    (3)            [below right= 10 mm  and 0 mm of 1]     {\small FBG strain data \\ \ary y}; 
	\node[pic node]   (4)     [above= 0mm of 1] {\includegraphics[width=\textwidth]{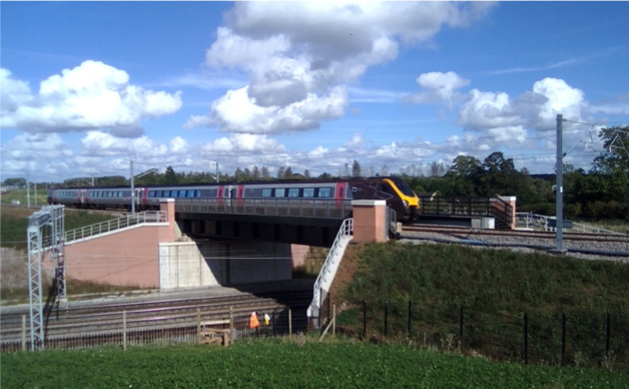}};
	\node[pic node]   (5)     [above= 0mm of 2] {\begin{minipage}[c]{\textwidth} \centering \includegraphics[width=0.9\textwidth]{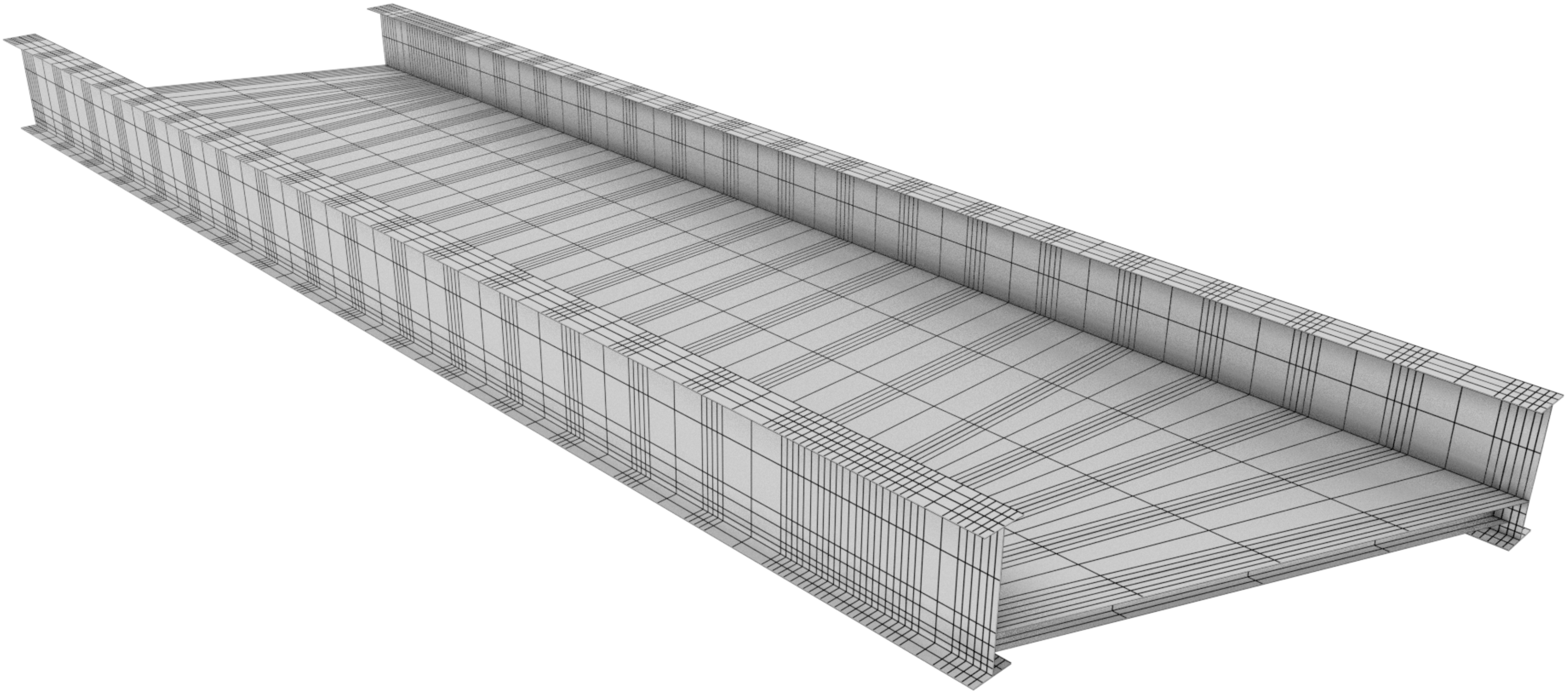} \end{minipage} };
	
	\path[every node/.style={font=\sffamily\small}] 
	(1) edge node [] {} (3)
	(3) edge node [] {} (2)
	(1) edge node [] {} (2) ;
	\end{tikzpicture} 
	\caption{Railway intersection bridge, a representative finite element discretisation and the statistical model underlying \mbox{statFEM}. The unobserved `true' bridge response~$\ary z$, the strain~$\ary y$ measured using FBG sensors and the finite element response~$\ary u$ are all random\label{fig:introFig}}
\end{figure*}

The statistical construction of the FE method, dubbed \mbox{statFEM},  recently introduced in~\cite{girolami2020statistical}  allows predictions to be made about the true system behaviour in light of limited sensor data and a misspecified FE model. As in all Bayesian approaches, any lack of knowledge, or uncertainty, associated with the errors in the collected data, choice of the FE  model and its parameters are represented as random variables. See~\cite{beck2010bayesian} or~\cite{huang2019state} for  an overview on Bayesian approaches in structural health monitoring. Typically, the as-built dimensions, support conditions and loadings of a structure are only partially known and can be treated as random.  Starting from some assumed (subjective) prior probability densities for the random variables, Bayes rule provides a coherent formalism to infer, or learn, the respective posterior densities using the likelihood of the observations~\citep{gelman2013bayesian, stuart2010inverse,kaipio2006statistical}. 

In \mbox{statFEM}, following~\cite{kennedy2001bayesian},  the observed data~$\ary y$ is equated to three random components: a FE component~$\ary u$, a model misspecification component~$\ary d$ and measurement noise~$\ary e$. The prior densities of each of the three random components may depend on additional random {\em hyperparameters} to be learned from the sensor data. In this paper, only the model misspecification component has hyperparameters. It is straightforward to introduce other hyperparameters pertaining to aspects of the FE model and the measurement noise, see~\cite{girolami2020statistical}. The choice of the additional hyperparameters leads to an increase in computational cost and requires particular care to circumvent non-identifiability~\citep{arendt2012quantification,nagel2016unified}.  The posterior density of the  FE component, the true system response and the hyperparameters are inferred via the Bayes rule. The overall approach is akin to the empirical Bayes or evidence approximation techniques prevalent in machine learning~\citep{mackay1992bayesian,mackay1999comparison,murphy2012machine}. The prior probability density of the FE component is obtained by solving a traditional probabilistic forward problem with random parameters~\citep{sudret2000stochastic,ghanem1991stochastic}.

This paper evaluates the application of \mbox{statFEM} to the development of a statistical digital twin of the superstructure of an operational railway bridge, see Figure~\ref{fig:introFig}.  The newly constructed skewed half-through bridge has previously been instrumented with fibre optic based sensors. The instrumentation of the bridge using FBG strain sensors is described in~\cite{butler2018monitoring}. The sensor measurements have since been compared with a deterministic FE model of the bridge in~\cite{lin2019performance}.  The contributions of the present work are threefold: (1) the first-time application of a new modelling paradigm for use in structural health monitoring;  (2) evaluating the trade-off between the number of sensor measurements and  the accuracy of statFEM prediction; and (3) introducing the concept of a statistical digital twin by demonstrating the application of \mbox{statFEM} to continuous strain sensing.

The outline of this paper is as follows. In Section~\ref{sec:selfSensingBridge} the structural system of the railway bridge and its instrumentation with FBG sensors are introduced. The proposed statistical digital twin of the superstructure is discussed in Section~\ref{sec:digitalTwin}. To this end, first the statistical finite element method is reviewed in Section~\ref{sec:statfem}, and the FE model of the structure is then  presented in Section~\ref{sec:problemDesc}. In Section~\ref{sec:results} the obtained digital twin is used to infer the true strain distribution over the entire structure during the passage of a train. In particular, it is shown that the inference results are largely insensitive to the number of sensors and sampling frequency confirming the utility of a \mbox{statFEM} model as an integral part of a digital twin.   Section~\ref{sec:conclusions} concludes the paper and discusses several promising directions for future research.

	%
\section{Self-sensing railway bridge} \label{sec:selfSensingBridge}

This section provides a summary of the structural system of the railway bridge and its sensor instrumentation with a fibre optic based sensor (FOS) network.
Both were previously introduced in~\cite{butler2018monitoring} and~\cite{lin2019performance}.

\subsection{Structural system } \label{sec:structuralSys}
\begin{figure}[]
	\centering
	\includegraphics[width=0.7\linewidth]{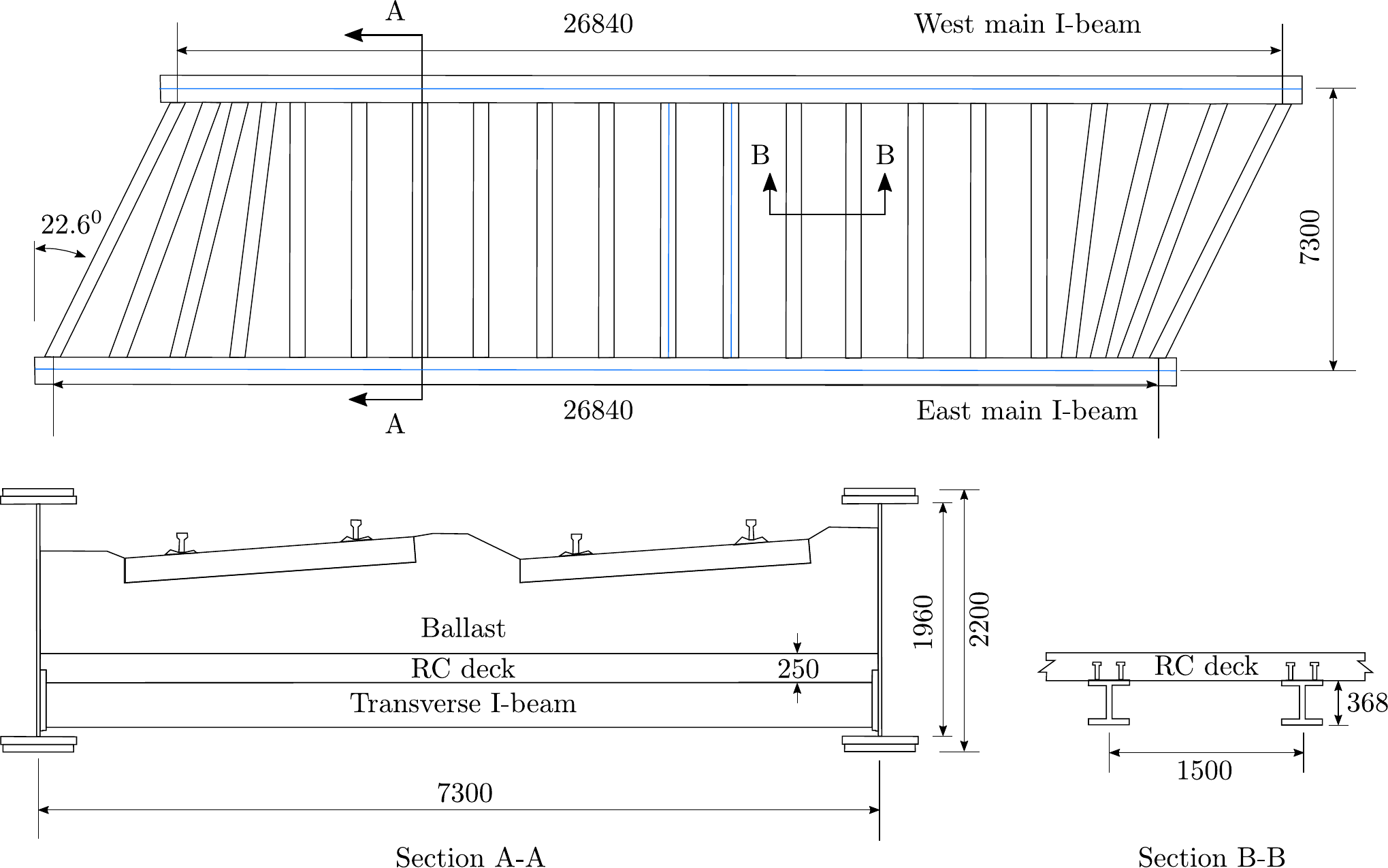}	
	\caption{Dimensions of the bridge superstructure (all in mm)}
	\label{fig:engDrawing}
\end{figure}

%
Completed in 2016, Intersection Bridge 20A is a $27.34~\text{m}$ steel skewed half-through railway bridge located along the West Coast Mainline in Staffordshire in UK (Figure~\ref{fig:introFig}).  The bridge carries two lines of passenger trains (with speeds up to $160~\text{km/h}$) over another heavily trafficked rail corridor along the West Coast Main Line. The superstructure consists of a pair of steel main I-beams and 21  transverse I-beams spanning the $7.3~\text{m}$ width of the structure, see Figure~\ref{fig:engDrawing}. The two main I-beams (referred to as east and west) are $26.84~\text{m}$ long and $2.2~\text{m}$ deep (including doubler plates). Web stiffeners welded along the outside web of the main I-beams are used to improve stability and to prevent local buckling of the web and top flanges. Four rocker-type bearings, which sit atop reinforced concrete abutments support the bridge superstructure. 

The transverse I-beams are $368~\text{mm}$ deep and the attached shear stud connectors provide for composite action with the reinforced concrete deck slab. The transverse beams are spaced at every $1.5~\text{m}$ in the middle of the bridge, and are fanned closer towards the two ends. Two  types of connection are used between the transverse and main I-beams. The transverse beams are consecutively either pinned with a 6-bolt end plate or moment connected with a 10-bolt stiffened end plate. A $250~\text{mm}$ thick reinforced concrete deck slab spans between transverse beams and supports the ballasted track bed system. The ballast has the minimum depth of $300~\text{mm}$ and supports the prestressed concrete sleepers to which the rails are fastened. 

%
\subsection{Sensor instrumentation \label{sec:sensorInstrumentation}}
\begin{figure*}[]
	\centering
	\includegraphics[width=0.8\linewidth]{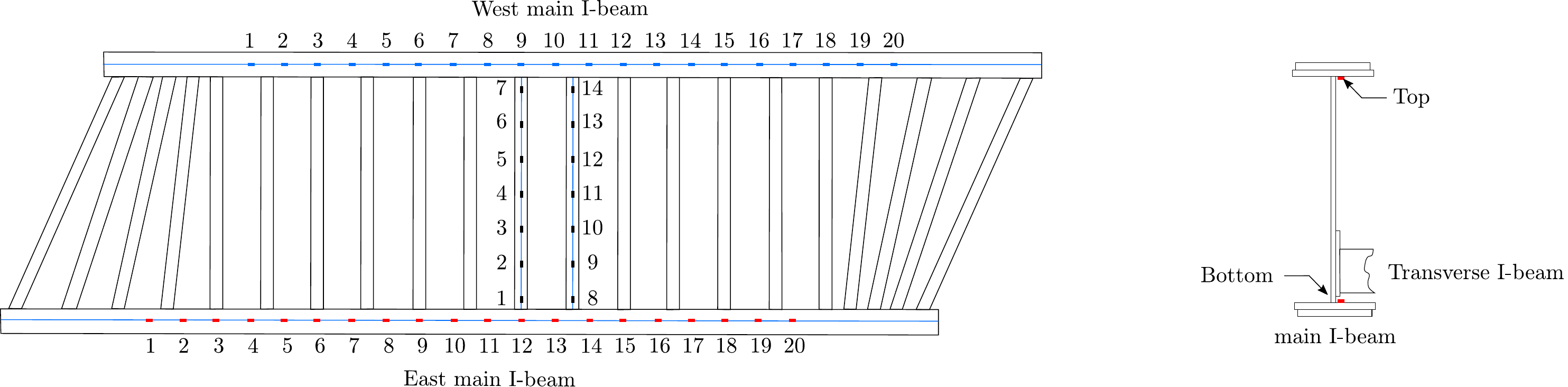}
	\caption{Position and numbering of the FBG strain sensors installed on the bridge. The same numbering is used for the sensors at  both top and bottom flanges. Not to scale}
	\label{fig:bridgeschematic}
\end{figure*}

%
This recently constructed bridge was instrumented with an innovative FOS network which was designed to provide reliable long-term measurements of the bridge's operational performance~\citep{butler2018monitoring}. The FBG sensors used were based on draw tower grating technology in which the gratings are inscribed at the time of glass drawing and prior to application of the core cladding. Compared with traditionally manufactured FBGs in which the core cladding must be stripped prior to FBG inscription, draw tower gratings are inherently more robust and, in the case of the sensors used in this project, enable the inscription of up to 20 FBGs along a single sensor cable. As depicted in Figure~\ref{fig:bridgeschematic}, FBG sensors were installed on the main structural elements. The sensor arrays used had up to 108 strain FBGs. Along the east main I-beam 20 FBGs were placed at one metre spacing along both the top and bottom flanges. Similarly, the same sensor arrangement is installed along the west main I-beam (i.e., 40 FBGs per I-beam). The other 28 FBGs are located along the top and bottom flanges of two adjacent transverse I-beams close to the midspan of the main I-beams. 

The sensors installed were based on FBG technology which can measure strain and temperature at discrete points along an optical fibre. Consisting of glass fibres, FOS are lightweight, require minimal number of wiring cables and  provide stable long-term measurements. When laser light is shone down a fibre optic cable, the gratings act like dielectric mirrors reflecting only those wavelengths, which match the Bragg wavelength. When the fibre optic cable is elongated (or shortened), the wavelength of the reflected light shifts in proportion to this change in strain. Once a fibre optic cable is attached to a structure, the cable can be used to measure changes in strain of the structure itself. The mechanical strain can be calculated by removing temperature effects associated with the thermal expansion of glass, the effect of temperature on the index of refraction, and the thermal expansion of the substrate material (i.e. steel) using
\begin{equation}
\Delta \epsilon_m = \frac{1}{k_\epsilon} \left(\left(\frac{\Delta \lambda}{\lambda_0}\right)_{S}-k_T \frac{\left(\frac{\Delta \lambda}{\lambda_0}\right)_T}{k_{T_T}}\right)-\alpha_{\text{sub}}\frac{\left(\frac{\Delta \lambda}{\lambda_0}\right)_T}{k_{T_T}} \, ,
\label{equation:strain}
\end{equation}
where $\Delta \epsilon_m$ is the change in mechanical strain; $(\Delta\lambda/ \lambda_0)_S$ is the change in relative wavelength of the strain sensor; $(\Delta\lambda/ \lambda_0)_T$ is the change in relative wavelength of the temperature-compensating sensor; $k_\epsilon$ is the gauge factor provided by the strain FBG manufacturer (typically 0.78); $k_T$ is  the change of the refractive index of glass; $k_{T_T}$ is the  experimentally derived constant for the FBG temperature-compensating sensor; and $\alpha_{\text{sub}}$ is the linear coefficient of thermal expansion of the substrate material \mbox{(concrete: 10$ \times 10^{-6}  / {^\circ\text{C}}$;} \mbox{steel: 12$ \times 10^{-6} / {^\circ\text{C}}$)}. Each FBG is capable of recording strain with an accuracy of approximately $\pm ~5 \cdot 10^{-6} $ at a data acquisition rate of up to $250~\text{Hz}$. 

Monitoring data were captured by the sensor network beginning from the time the bridge was being constructed~\citep{butler2018monitoring}. Following the construction, the bridge was monitored sporadically for two years~\citep{lin2019performance}. During this period, FBG data were captured during the passage of more than 130 individual trains. This dataset contains strain data for the passage of two passenger train types: the London Midland Class 350 Desiro (type T1) and the Cross Country Class 221 Super Voyager (type T2). Both train types included 4 or 5-car configurations. Their axle weightings and spacing were provided by Network Rail, the UK's national rail authority.

	%
\section{Statistical Digital Twin \label{sec:digitalTwin}}
%
This section provides a brief review of statFEM and presents the FE model of the bridge superstructure. For further details on statFEM and its application to diffusion problems refer to~\cite{girolami2020statistical}. In the FE model the main sources of uncertainty are the partially known train weight, the resulting structural loads and the particulars of the FE model.  The material properties and geometry parameters are assumed to be deterministic. Moreover, the inertia effects are not taken into account due to the relatively short span of the bridge~\citep{lin2019performance}.

\subsection{Review of the statistical finite element method  \label{sec:statfem}}
%

%
\subsubsection{Probabilistic forward FE formulation  \label{sec:forwardModel}}
%
The bridge superstructure is modeled using isogeometrically discretised Kirchhoff-Love shell finite elements~\citep{Cirak:2000aa,Cirak:2002aa,Cirak:2011aa}. The shell model takes into account the in-plane membrane and out-of-plane bending response of the structural members. The FE model of the superstructure, consisting of the main girders, cross I-beams and the reinforced concrete deck, is obtained by rigidly connecting horizontally and vertically aligned plates, i.e. initially planar shells, along joints. The joints are able to transfer both forces and moments. 

The weak form of the shell equilibrium equation, or the principle of virtual work, for a Kirchhoff-Love shell with the midsurface~$\Omega$ and the position vector~$\vec x$ reads
\begin{equation} \label{eq:weak}
	\int_{\Omega}  \left ( \vec n (\vec x) : \delta \vec \alpha  (\vec x) + \vec m (\vec x) : \delta \vec \beta (\vec x) \right )  \D \Omega = \sum_j \vec f^{(j)} \cdot  \delta \vec u^{(j)} + \int_\Omega   \vec r (\vec x)  \cdot \delta \vec u (\vec x)  \D \Omega  \, . 
\end{equation}
Here the penalty term enforcing the conformity of the displacements and rotations of all plates attached to the same joint has been omitted for the sake of brevity.  The internal virtual work on the left-hand side depends on the membrane force resultant~$\vec n(\vec x)$, the bending moment resultant~$\vec m(\vec x)$, the virtual membrane strain~$\delta \vec \alpha(\vec x)$ and the virtual bending strain~$\delta \vec \beta(\vec x)$. Evidently, these fields depend in turn either on the displacements~$\vec u(\vec x)$ or the virtual displacements~$\delta \vec u(\vec x)$. Furthermore, the material is assumed to be linear elastic and isotropic so that the two resultants~$\vec n(\vec x)$ and~$\vec m(\vec x)$ depend on the respective strains~$\vec \alpha(\vec x)$ and~$\vec \beta(\vec x)$ via the Young's modulus~$E$, Poisson's ratio~$\nu$ and the shell thickness~$h$. The right-hand side of~\eqref{eq:weak} represents the external virtual work and depends on the deterministic concentrated forces~$\vec f^{(j)}$ applied at the respective positions~$\vec x^{(j)}$ and the random distributed loading~$\vec r (\vec x)$.  The random distributed loading takes into account the uncertainty in train weight and other uncertainties pertaining to the choice of the finite element model.  As will be detailed later, it is assumed that the loading applied by the train axles is composed only of vertical components so that both~$\vec f^{(j)}$ and~$\vec r(\vec x)$ have only a non-zero vertical component. The random distributed train loading~$\vec r(\vec x)$ is a Gaussian process with a zero mean and prescribed covariance, i.e.  
\begin{equation} \label{eq:randomR}
	\vec r (\vec x) = 
	\begin{pmatrix} 
	0 \\
	0 \\
	\set {GP} \left ( 0, \, c_r (\vec x, \, \vec x') \right ) 
	\end{pmatrix}  \, . 
\end{equation}
Without loss of generality, the squared-exponential kernel is chosen as the covariance function
\begin{align} \label{eq:forcingKernel}
	c_r(\vec x, \vec x') = \sigma_r^2 \exp \left(  - \frac{ \| \vec x - \vec x' \|^2 }{2 \ell_r^2}   \right) \, , 
\end{align}
where~$\sigma_r$ is a scaling factor and~$\ell_r$ is a lengthscale parameter. Although the Euclidean squared-distance~\mbox{$\| \vec x - \vec x' \|^2$} has been used, in structural mechanics it may be more appropriate to use the squared geodesic distance~\citep{scarth2019random}. 

The position and displacement vectors~$\vec x$ and~$\vec u (\vec x)$ in the weak form~\eqref{eq:weak} are discretised using basis functions obtained from Catmull-Clark subdivision surfaces~\citep{Cirak:2000aa,zhang2018subdivision}. The weak form~\eqref{eq:weak} depends on the curvature so that the FE basis functions must be smooth or, in other words, must have square-integrable second derivatives. Subdivision surfaces are the generalisation of B-splines from computer-aided geometric design (CAGD) to unstructured meshes and provide smooth basis functions for discretising the position and displacement vectors~$\vec x$ and~$\vec u(\vec x)$. The resulting FE approach is referred to as isogeometric analysis~\citep{Hughes:2005aa} and allows one to use the same representation for geometric design and analysis. Without going into further detail,  the subdivision basis functions~$\phi_i(\theta^1, \, \theta^2)$ furnish the two approximants
\begin{align}  \label{eq:approx}
	\vec x(\theta^1, \, \theta^2) = \sum_i \phi_i (\theta^1, \, \theta^2) \vec x_i  \, , \qquad   \vec u (\theta^1, \, \theta^2) = \sum_i  \phi_i (\theta^1, \, \theta^2) \vec u_i \, ,
\end{align}
where~$\vec x_i$ is the coordinate and~$\vec u_i$ is the displacement of the FE node with the index~$i$ and $(\theta^1, \, \theta^2)$ are the parametric coordinates. The sums in~\eqref{eq:approx} are over all the nodes in the FE mesh. After introducing~\eqref{eq:approx} into the weak form~\eqref{eq:weak} and numerically evaluating the integrals the linear system of equations 
\begin{equation}
	\ary A \ary u = \ary f
\end{equation}
is obtained. Herein,~$\ary A$ is the stiffness matrix, $\ary u$ is the nodal displacement vector and the force vector~$\ary f$ has  the multivariate Gaussian density 
\begin{align} \label{eq:randomF}
	\ary f \sim p(\ary f) = \set N \left(\overline{\ary f}, \, \ary C_{\ary f} \right)  \, ,
\end{align}
where~$\overline{\ary f}$ is the mean vector and $\ary C_{\ary f}$ is the covariance matrix. See~\cite{girolami2020statistical} for the computation of~$\overline{\ary f}$ and~$\ary C_{\ary f}$. Because the stiffness matrix~$\ary A$ is deterministic  it is easy to show that the resulting random displacement vector has the probability density 
\begin{align} \label{eq:randomU}
\ary u \sim p(\ary u) = \set N \left(\ary A^{-1} \overline{\ary f}, \, \ary A^{-1} \ary C_{\ary f} \ary A^{-\trans} \right) 
                                =  \set N \left(\overline{\ary u}, \,  \ary C_{\ary u} \right) \, .
\end{align}
Although not considered in the present study, it is possible to consider uncertainties in the internal work in~\eqref{eq:weak}, such as random material parameters or geometry, leading to a random stiffness matrix~$\ary A$.  

%
\subsubsection{Statistical model and Bayesian inference  \label{sec:bayesianInference}}
%
In the posited statistical model, the observed strain $\ary y$ at the~$n_y$ sensor locations is assumed to be composed of three random components, i.e.,
\begin{align}\label{eq:statModel}
	\ary y = \ary z + \ary e = \rho \ary P \ary u + \ary d + \ary e   = \ary z + \ary e \, .
\end{align}
The observed strain~$\ary y$ is equal to the  sum of the unknown (i.e. unobserved) `true' strain~$\ary z$ and the Gaussian measurement error 
\begin{equation}
	\ary e  \sim p(\ary e) = \set N \left(\ary 0, \, \ary C_{\ary e} \right)
\end{equation}
with a diagonal covariance matrix~$\ary C_{\ary e} = \sigma_e^2 \ary I$ and a standard deviation~$\sigma_e$. In turn, the true system response~$\ary z$  is the linear combination of the FE strain~$\ary P \ary u$, depending on the nodal displacements~$\ary u$, and the mismatch, or model inadequacy, error~$\ary d$. In the FE component,~$\rho$ is an unknown scaling parameter and~$\ary P$ is a matrix for obtaining the strain at the~$n_y$ sensor locations from the nodal displacement vector~$\ary u$. The three random components~$\ary u$, ~$\ary d$ and~$\ary e$ are assumed to be statistically independent.

The random parameter~$\rho$ and the random mismatch error~$\ary d$ are the unknowns of the statistical model~\eqref{eq:statModel} and will be characterised using the measured strain data~$\ary y$ and the random FE solution~$\ary u$ in~\eqref{eq:randomU}. The probability density of the mismatch error is assumed to be a multivariate Gaussian 
\begin{equation}
	\ary d \sim p (\ary d) = \set N \left(\ary 0, \, \ary C_{\ary d} \right) \, ,
\end{equation}
and the covariance matrix~$\ary C_{\ary d}$ is obtained from the squared-exponential kernel 
\begin{equation} \label{eq:mismatchCov}
	c_d(\vec x, \vec x') = \sigma_d^2 \exp \left(  - \frac{ \| \vec x - \vec x' \|^2 }{2 \ell_d^2}   \right) \, 
\end{equation}
with the hyperparameters~$\sigma_d$ and~$\ell_d$. In the following, the three hyperparameters of the statistical model are collected in the vector 
\begin{equation}
	\ary w  \coloneqq \begin{pmatrix} \rho & \ell_d & \sigma_d \end{pmatrix}^\trans \, . 
\end{equation}

Next, the Bayesian inference of the random finite element solution~$\ary u$ and the hyperparameters~$\ary w$ in light of the observed data~$\ary y$ and the statistical model~\eqref{eq:statModel} are considered. Since all the variables in the statistical model are Gaussians one can write by inspection for the likelihood, see e.g.~\cite{murphy2012machine}, 
\begin{equation} \label{eq:likelihood}
	p(\ary y| \ary u, \, \ary w ) = \set N \left(\rho \ary P \ary u, \ary C_{\ary d} + \ary C_{\ary e} \right)  \, .
\end{equation}
According to the Bayes rule, the posterior density of the FE solution is given by
%
\begin{align}
	p( \ary u | \ary y , \, \ary w) &= \frac{p(\ary y| \ary u, \, \ary w ) p(\ary u ) }{p (\ary y | \ary w)}  \label{eq:postU1}
\end{align}
%
Note that in the present study the prior~$p(\ary u)$ does not have any hyperparameters, i.e. \mbox{$ p(\ary u, \,  \ary w) =  p(\ary u)$}. 

In statFEM one is interested in the posterior FE density 
\begin{equation}
	p (\ary u | \ary y) = \int p(\ary u  | \ary w , \, \ary y) p(\ary w | \ary y)  \D \ary w \, .
\end{equation}	
Following an empirical Bayes or evidence approximation approach, see e.g.~\cite{mackay1999comparison},  this integral can be approximated by 
\begin{align} 
	p(\ary u | \ary y) \approx p (\ary u |  \ary w^* , \,  \ary y)   \, . 
 \end{align}
As a point estimate~$\ary w^*$ either the maximum 
\begin{subequations}  \label{eq:pointEstimate}
\begin{align} 
	\ary w^* = \argmax_{\ary w} p(\ary y | \ary w)
	\intertext{or the mean}
	 \ary w^* = \expect [ p(\ary y | \ary w) ]  
\end{align}
\end{subequations}	 
are suitable. When the hyperparameters are endowed with a prior~$p(\ary w)$  the posterior \mbox{$p(\ary w | \ary y) \propto p(\ary y| \ary w) p(\ary w)$} can be considered for obtaining the point estimate. In this paper~$p(\ary w | \ary y)$ is sampled with the Markov-chain Monte Carlo (MCMC) method and the empirical mean of the samples is taken as~$\ary w^* $. According to the statistical model~\eqref{eq:statModel} and noting that linear combinations of Gaussian vectors is a Gaussian the marginal likelihood is given by
\begin{equation}
	p(\ary y| \ary w) = \set N \left(\rho \ary P \overline{\ary u}, \, \ary C_{\ary d} + \ary C_{\ary e} + \rho^2 \ary P \ary C_{\ary u}  \ary P^\trans \right)  \, .
\end{equation}

Finally, introducing the obtained point estimate~$\ary w^*$ in~\eqref{eq:postU1} yields the sought posterior density~$p(\ary u | \ary y)$. To this end, note that the densities on the right-hand side of~\eqref{eq:postU1} are Gaussians so that the posterior is a Gaussian as well. As shown in~\cite{girolami2020statistical}, the posterior density reads
\begin{subequations}\label{eq:posteriorU}
\begin{align}
	p(\ary u | \ary y ) &= \set N \left(\overline{\ary u}_{\ary y}, \, \ary C_{\ary u|\ary y} \right)
	\intertext{with the covariance matrix}
	\ary C_{\ary u| \ary y} &=  \ary C_{\ary u} - \ary C_{\ary u} \ary P^\trans \left ( \frac{1}{\rho^2} \left ( \ary C_{\ary d} + \ary C_{\ary e} \right )  +  \ary P \ary C_{\ary u} \ary P^\trans \right )^{-1} \ary P    \ary C_{\ary u}  
	\intertext{and the mean}
	\overline {\ary {u}}_{| \ary y} &=   \ary C_{\ary u| \ary y}   \left(  \rho  \ary P^\trans  \left ( \ary C_{\ary d} + \ary C_{\ary e} \right )^{-1} \ary y   +   \ary C_{\ary u}^{-1}  \overline{\ary u}  \right)  \label{eq:posteriorUc} \, .
\end{align}
\end{subequations}

The mean can also be expressed as
\begin{equation}
	\overline {\ary {u}}_{| \ary y} = \overline{\ary u} +  \ary C_{\ary u} \ary P^\trans \left ( \frac{1}{\rho^2} \left (  \ary C_{\ary d} + \ary C_{\ary e} \right )  + \ary P \ary C_{\ary u} \ary P^\trans \right )^{-1} \left ( \frac{\ary y}{\rho} -  \ary P \overline{ \ary u} \right ) \, .
\end{equation}
That is, the posterior mean is obtained by updating the prior mean~$\overline{\ary u}$ with the difference between the observed and prior mean strain~$\ary y / \rho - \ary P \overline{\ary u} $. The respective weight depends on the scaling parameter~$\rho$ and the relative magnitude of the covariance matrices~$\ary C_{\ary u}$, $\ary C_{\ary d}$ and~$\ary C_{\ary e}$.

With the determined posterior FE density the posterior true strain density is given by
\begin{align}\label{eq:posteriorZ}
	p(\ary z | \ary y) = \set N \left( \rho \ary P \overline{\ary u}_{| \ary y}, \rho^2 \ary P \ary C_{\ary u | \ary y} \ary P^{\trans} + \ary C_{\ary d} \right) \, .
\end{align}

%
\subsubsection{Mismatch modelling and multiple observations \label{sec:discrepancyAndMobserv}}
%
The previous section considered only one single observation vector~$\ary y$ recorded at a fixed time instant~$t_k$. However, the passage of a single train yields several hundreds of such observation vectors. These vectors are collected in the observation matrix
\begin{align}
	\ary Y = \begin{pmatrix} \ary y_0 & \ary y_1 & \dotsc  &  \ary y_{n_o}  \end{pmatrix} \, .
\end{align}
In principle, it is possible to consider each time instant~$t_k$ independently and to obtain for each a point estimate~$\ary w^*$ according to~\eqref{eq:pointEstimate}. This can be very costly given that each time instant requires MCMC sampling. More critically, the information content available in one single observation vector~$\ary y$ often gives  only a very poor estimate for~$\ary w^*$. 
Consequently, it is advantageous to consider all observation vectors simultaneously and to obtain one single~$\ary w^*$. To this end, the mismatch error covariance kernel~\eqref{eq:mismatchCov} is replaced with the scaled squared-exponential kernel
\begin{subequations}
\begin{align} \label{eq:mismatchCovT}
	c_d(\vec x, \, \vec x', \, k) &= \left(\gamma_k  \sigma_d \right)^2 \exp \left(  - \frac{ \| \vec x - \vec x' \|^2 }{2 \ell_d^2}   \right) 
	\intertext{with}
	 \gamma_k &= \frac{\| \ary f_k \|}{ \max_k \| \ary f_k \| } \, ,
\end{align}
\end{subequations}
where the index~$k$ corresponds to the time instant~$t_k$ and $\| \cdot \|$ denotes~$L_2$ norm. The scale factor~$\gamma_k$ is necessary because the number of axles on the bridge, hence the loading~$\ary f_k$ at each~$t_k$ is different. As implied by the statistical model~\eqref{eq:statModel}, it is necessary to scale each random vector~$\ary d_k$ in dependence of the magnitude of the observed strain~$\ary y_k$, which is approximately proportional to the actual loading~$\ary f_k$. Assuming statistical independence between the observations, the marginal likelihood for determining~$\ary w^*$ reads 
\begin{subequations} \label{eq:marginalY}
\begin{align}
	p(\ary Y | \ary w) &= \prod_{k=1}^{n_o} p(\ary y_k | \ary w) 
	\intertext{with}
	p(\ary y_k| \ary w) &= \set N \left(\rho \ary P \overline{\ary u}_k, \, \ary C_{\ary d_k} + \ary C_{\ary e} + \rho^2 \ary P \ary C_{\ary u_k}  \ary P^\trans \right)  \, . 
\end{align}
\end{subequations}
%

The prior mean~$ \overline{\ary u}_k$ and covariance~$\ary C_{\ary u_k}$ are given by~\eqref{eq:randomU} and correspond to the loading~$\ary f_k$ at the time~$t_k$.

In summary, starting with the FE prior densities~$p (\ary u_k) = \set N ( \overline{\ary u}_k , \,  \ary C_{\ary u_k})$ and the observation matrix~$\ary Y$ the posterior FE and true strain densities are determined in three steps.
\begin{enumerate}
	\item Sample the marginal likelihood~$p(\ary Y | \ary w)$ given by~\eqref{eq:marginalY} using MCMC.
	\item Determine from the MCMC samples the point estimate~$\ary w^*$.
	\item Compute the posterior FE and true strain densities~$p(\ary u_k | \ary y_k)$ and~$p(\ary z_k | \ary y_k )$ by introducing the determined~$\ary w^*$ into \eqref{eq:posteriorU} and~\eqref{eq:posteriorZ} .   
\end{enumerate}

%
\subsection{FE model of the bridge superstructure \label{sec:problemDesc}}
%

\begin{figure}[]
	\centering
	\includegraphics[width=0.55\linewidth]{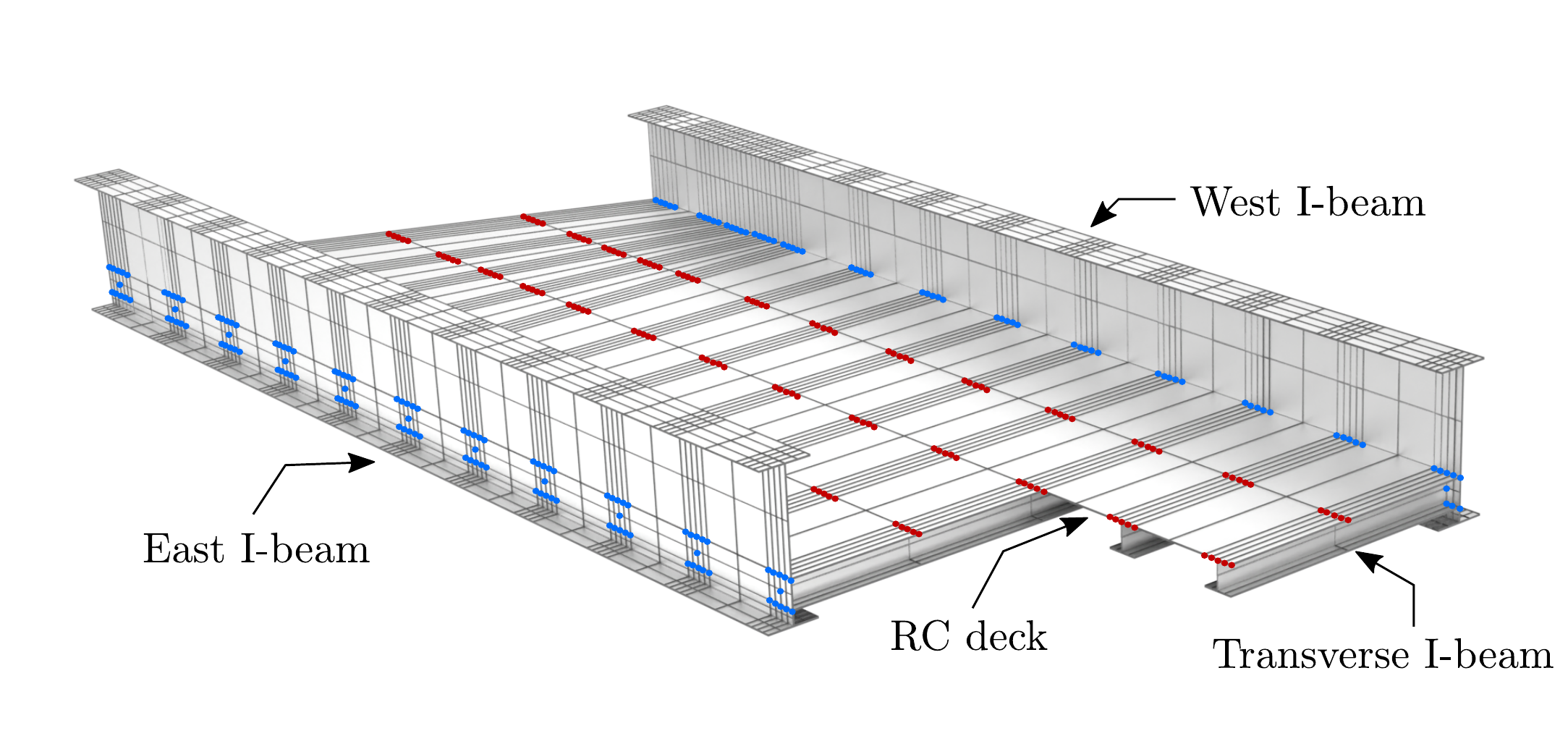}	
	\caption{Finite element model of the bridge superstructure. The two main I-beams and the 21 transverse I-beams are moment connected at the mesh nodes in blue. The reinforced concrete deck is moment connected to the transverse I-beams at the mesh nodes in red}
	\label{fig:constraints}
\end{figure}

As detailed in Section~\ref{sec:structuralSys}, the bridge superstructure is comprised of two identical main I-beams (east and west), 21 transverse I-beams, and a reinforced concrete deck. The transverse I-beams are modelled as rigidly connected to the reinforced concrete deck owing to the double row of shear connectors atop the transverse I-beams. In the actual structure the transverse beams are alternatingly pinned and moment connected to the main I-beams. For simplicity in the considered model, following the numerical studies in~\cite{lin2019performance}, all transverse beams are modelled as  moment connected.  Consequently, all components of the superstructure are modelled using shell finite elements and along joints moment connected to each other, see Figure~\ref{fig:constraints}. Each moment connected joint is able to transfer both forces and moments.  As noted in Figure~\ref{fig:engDrawing}, in the FE model the two main I-beams are $26.84~\text{m}$ long, $0.7~\text{m}$ wide and~$2.04~\text{m}$ deep. The two main I-beams are~$7.3~\text{m}$ apart.  Instead of doubler plates as in the actual structure, the flanges consist of a single plate with equivalent second moment of area.  The 21 transverse I-beams are categorised either as orthogonal or fanned. The orthogonal beams are located in the centre of the superstructure and are placed at every~$1.5~\text{m}$; and, the fanned beams are situated near the two ends of the bridge. The transverse beams are $0.4~\text{m}$ wide and $0.4~\text{m}$ deep. The dimensions of all the I-beams is listed in Table~\ref{tab:geometryComponents}. 

\begin{table}[b]
\begin{tabular}{lcccc}
\cline{2-5}
               & \multicolumn{2}{c}{main} & \multicolumn{2}{c}{transverse} \\
               & web        & flange      & web           & flange         \\ \midrule
thickness      & 25         & 120         & 16.5          & 27             \\
depth or width & 2040       & 700         & 400           & 400            \\ \bottomrule
\end{tabular}
			\caption{Web and flange dimensions of the main and transverse I-beams (all in $\text{mm}$). \label{tab:geometryComponents}}%
\end{table}

The Young's modulus and Poisson's ratio of the steel I-beams are assumed to be $E_s = 210~\text{GPa}$ and $\nu_s = 0.3$, respectively. The deck of the bridge is modelled without the ancillary structures such as the ballast and prestressed concrete sleepers. The thickness of the deck in the FE model is $0.25~\text{m}$. Assuming a uniform distribution of the steel reinforcement across the cross-section of the concrete deck the equivalent reinforced concrete elasticity modulus according to the rule of mixtures is given by 
\begin{align} \label{eq:rcModulus}
E_{rc} = q E_s + (1-q) E_c \, .
\end{align}
A steel reinforcement ratio of \mbox{$q = 0.03$} and a concrete elasticity modulus of \mbox{$E_c = 35~\text{GPa}$} are assumed yielding \mbox{$E_{rc} = 40.5~\text{GPa}$.} The Poisson's ratio of the reinforced concrete is \mbox{$\nu_{rc} = 0.2$.}   The equivalent elasticity modulus~\eqref{eq:rcModulus}  does not take into account that the rebar is usually placed away from the centre of the cross-section. 

\begin{figure} []
	\centering
	\subfloat[Top view]{
		\includegraphics[width=0.35\textwidth]{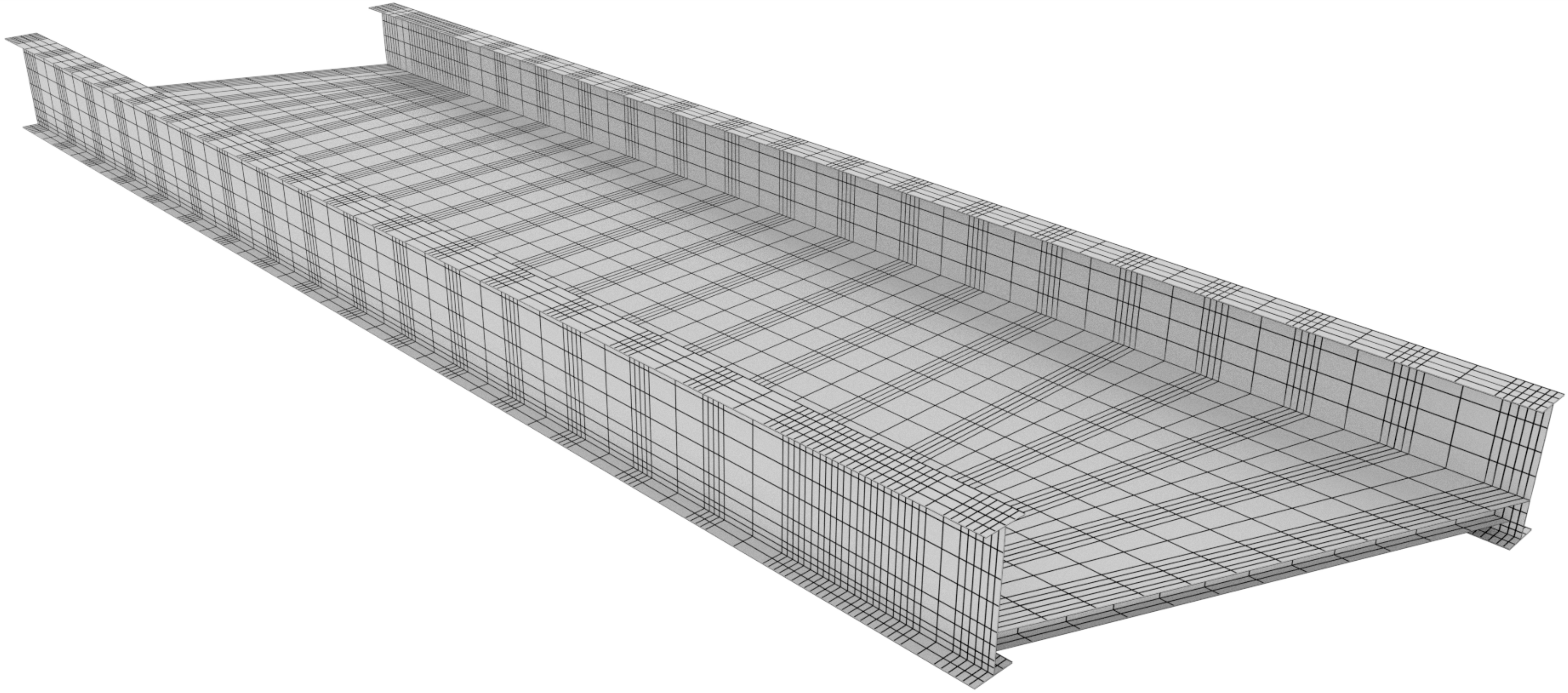} \label{fig:mesh2top}
	} 	\qquad
	\subfloat[Bottom view]{
	\includegraphics[width=0.3\textwidth]{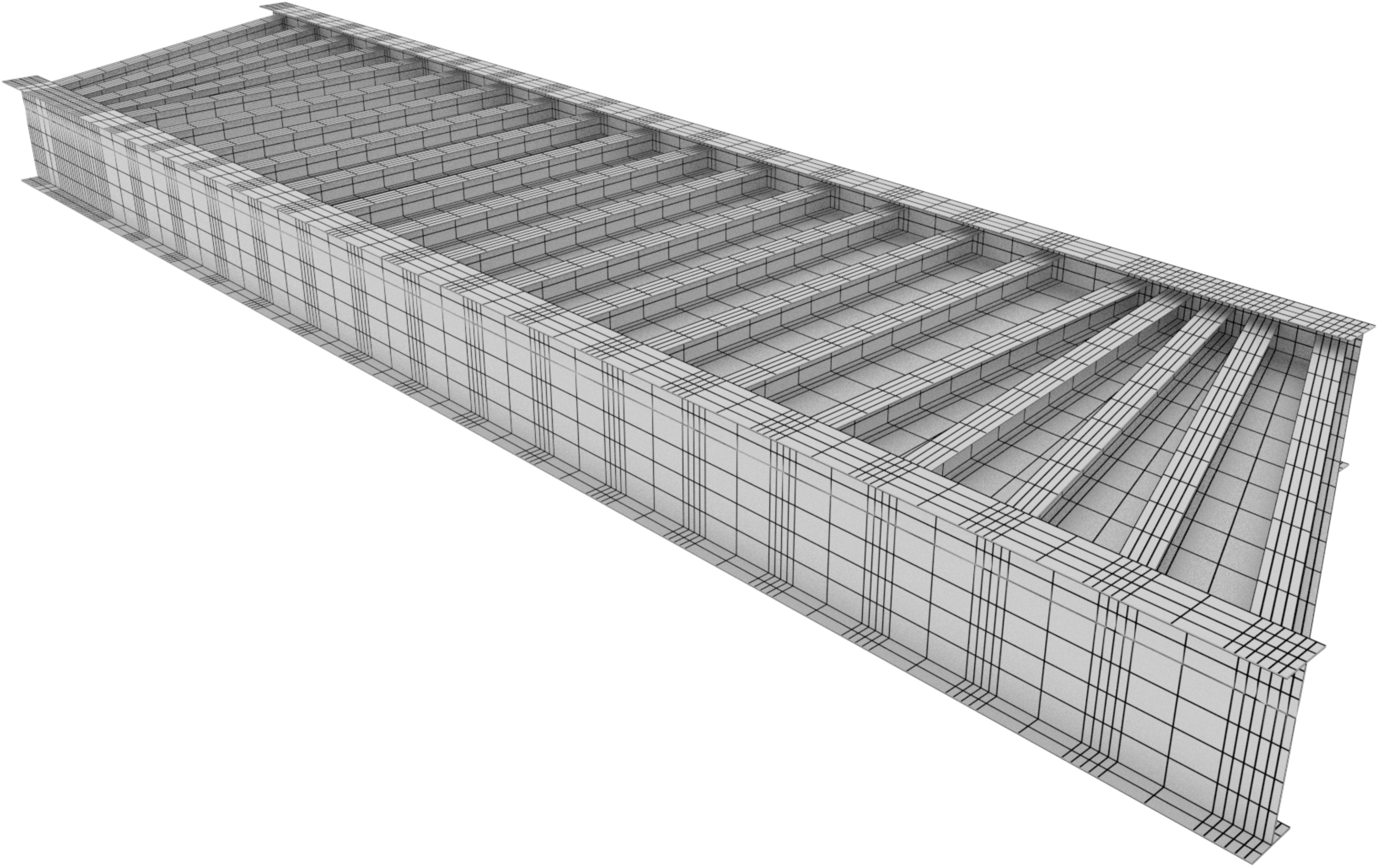} \label{fig:mesh2bottom}
	} 
	\caption{ The fine FE mesh M2} \label{fig:mesh2}
\end{figure}

Two different finite element meshes for discretising the bridge superstructure are considered. Mesh M1 is a relatively coarse mesh consisting of 4635 nodes and 4600 quadrilateral elements and the finer mesh M2 consisting of 7191 nodes and 7144 quadrilateral elements, see Figures \ref{fig:introFig} and \ref{fig:mesh2}. To model the rocker type bridge supports, pinned and roller supports at the relevant FE nodes are applied. 

\begin{figure*} []
	\setlength{\fboxsep}{0pt}%
	\setlength{\fboxrule}{0pt}%
	\centering
	\includegraphics[width=0.55\textwidth]{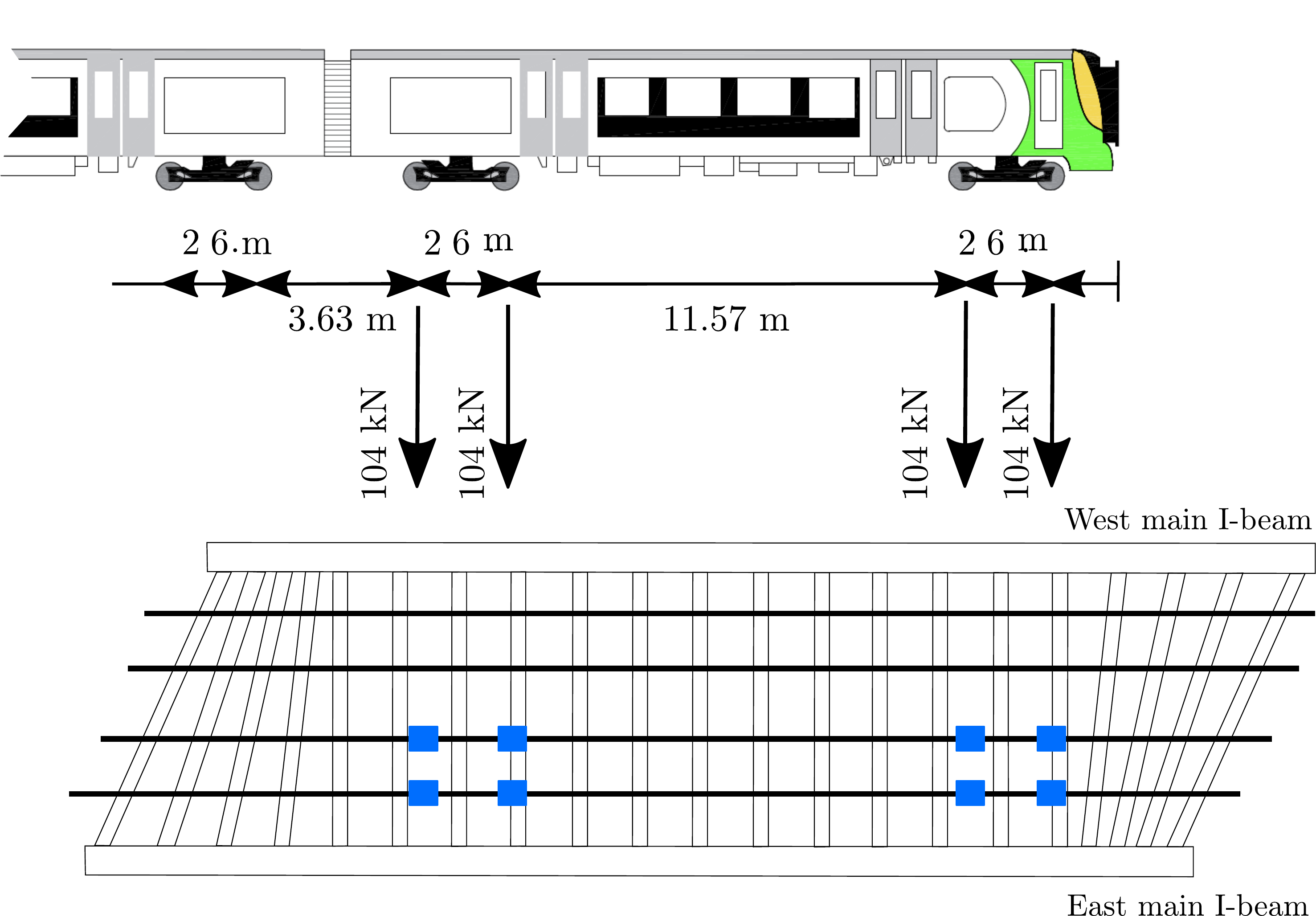} 
	\caption{Location and magnitude of the forces applied to the superstructure. Forces are directly applied to the deck. The axle weight of 104 kN is split into two forces representing the two attached wheels}  \label{fig:axleDef} 
\end{figure*}

Although there is data for two different types of train crossing the bridge, this paper focuses on type T1 train (London Midland Class 350 Desiro) with four carriages and a total length of $81.47~\text{m}$~\citep{lin2019performance}. The average normal axle load is $104~\text{kN}$ ($52~\text{kN}$ per wheel), see  Figure~\ref{fig:axleDef}. The speed of the train is approximately $131~\text{km/h}$ so that the train crosses the bridge in less than three seconds. The moving loads generated by the train wheels are directly applied to the concrete slab neglecting any load distribution effects by the sleepers and the track ballast.  Moreover, the rail tracks are assumed to have no vertical cant so that any centrifugal forces are neglected. Clearly, a more faithful modelling of the train loading and its transfer to the bridge superstructure requires a significantly more advanced FE model.   
Due to the relatively short span of the bridge, inertia effects are not taken into account in the FE analysis.  That is, for a given set of wheel forces and positions a single static analysis is performed. 

The weight of the train and its approximation as static point loads in the FE model is obviously subject to a number of assumptions and uncertainties. In the proposed Bayesian statFEM framework these uncertainties are taken into account by considering the FE solution~$\ary u$ as a random variable. It bears emphasis that this randomness is epistemic, i.e. it is due to lack of knowledge and not due to inherent variability of the FE solution.  To obtain the respective FE prior density~$p(\ary u)$ in~\eqref{eq:randomU}, a random loading~$\vec r (\vec x)$  as described by the Gaussian process~\eqref{eq:randomR} is applied. The respective covariance scaling and length scale parameters are chosen as \mbox{$\sigma_r = 1000~\text{Pa}$} and $\ell_r = 1~\text{m}$. 
The scaling factor $\sigma_r$  represents  20\% of  the maximum train load distributed over the planform of the bridge, see Figure~\ref{fig:forceOverTime}. The length scale parameter $\ell_r$ is chosen to be in the order of the  distance between the transverse beams and is about one-third of the ballast depth.  The random loading~$\vec r(\vec x)$  is applied in addition to the deterministic point forces of~$52~\text{kN}$ per wheel. The choice of the two parameters~$\sigma_r$ and~$\ell_r$  requires input from domain experts, i.e. bridge engineers, and may be formalised with prior elicitation techniques from statistics, see the recent review~\cite{mikkola2021prior}. Alternatively,~$\sigma_r$ and~$\ell_r$ can be interpreted as hyperparameters and inferred from the observed data~\citep{kennedy2001bayesian,nagel2016unified}.

\section{Results and discussion \label{sec:results}}
\begin{figure}[]
	\centering
	\subfloat[Total train loading over time]{
		\includegraphics[width=0.42\textwidth]{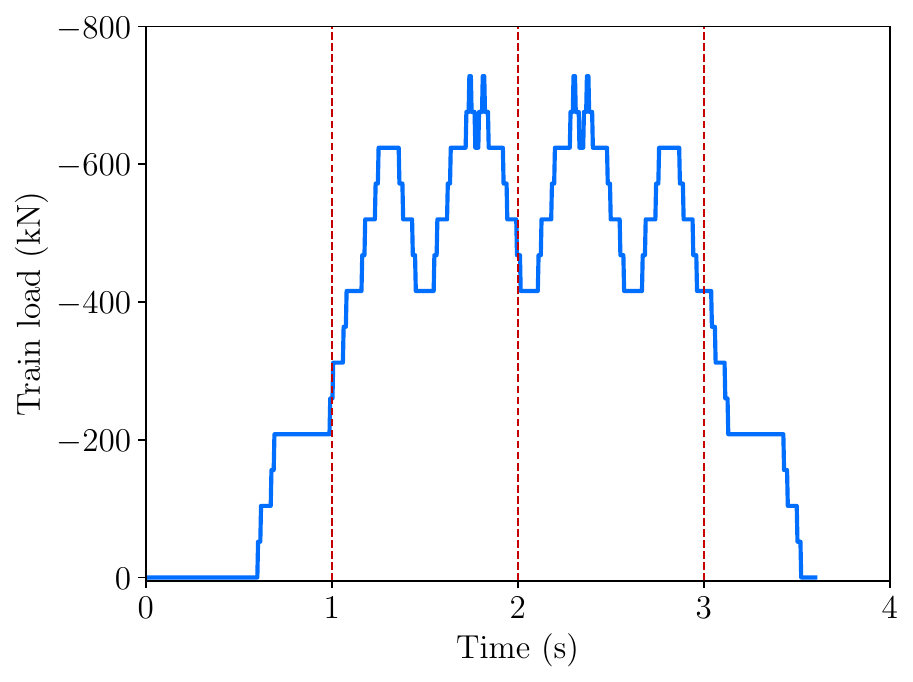} \label{fig:forceOverTime}
	} \\
	\subfloat[$t = 1 \, s$]{
		\includegraphics[width=0.3\textwidth]{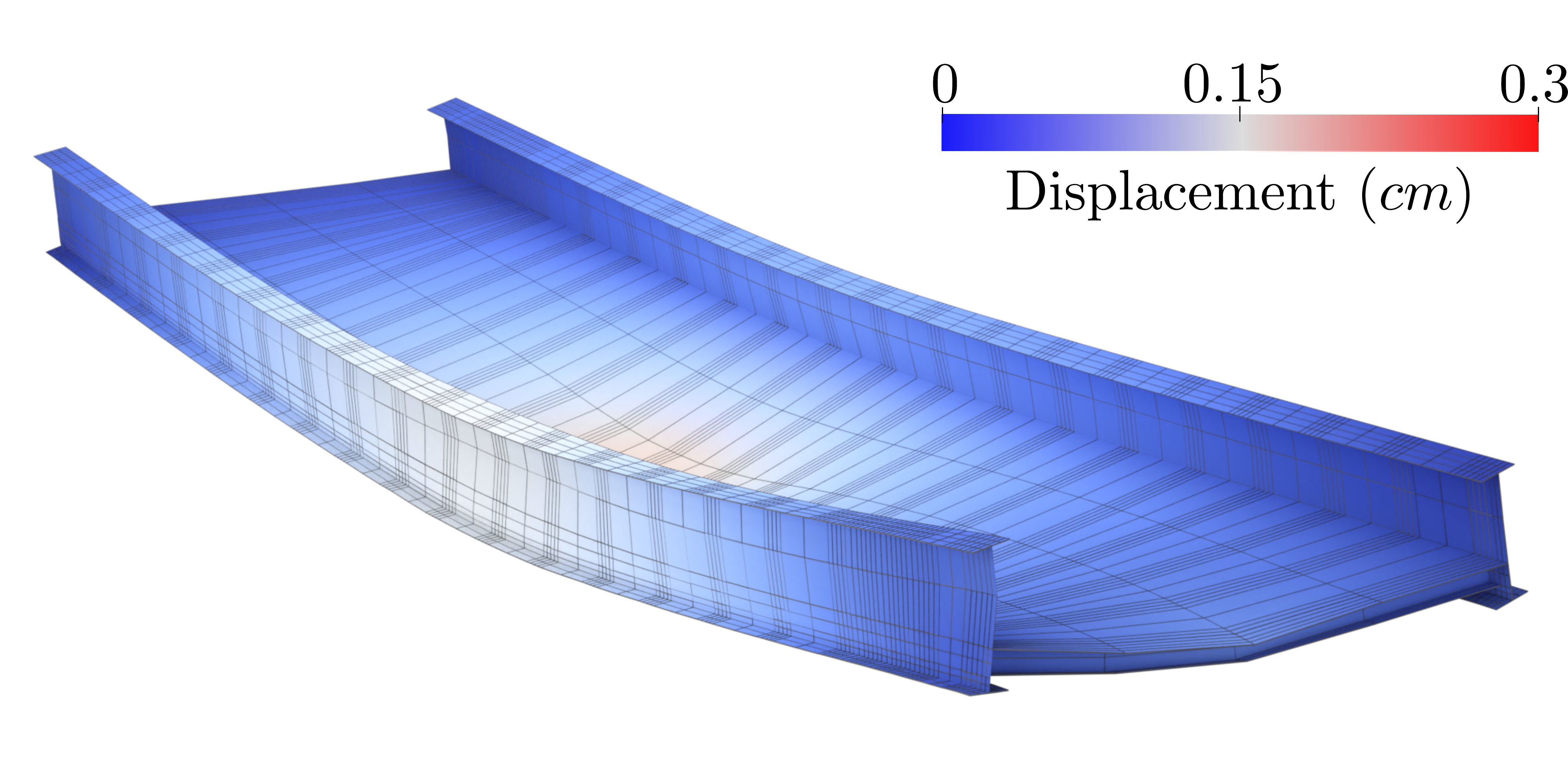} 
	} 
	\subfloat[$t = 2 \, s$]{
		\includegraphics[width=0.3\textwidth]{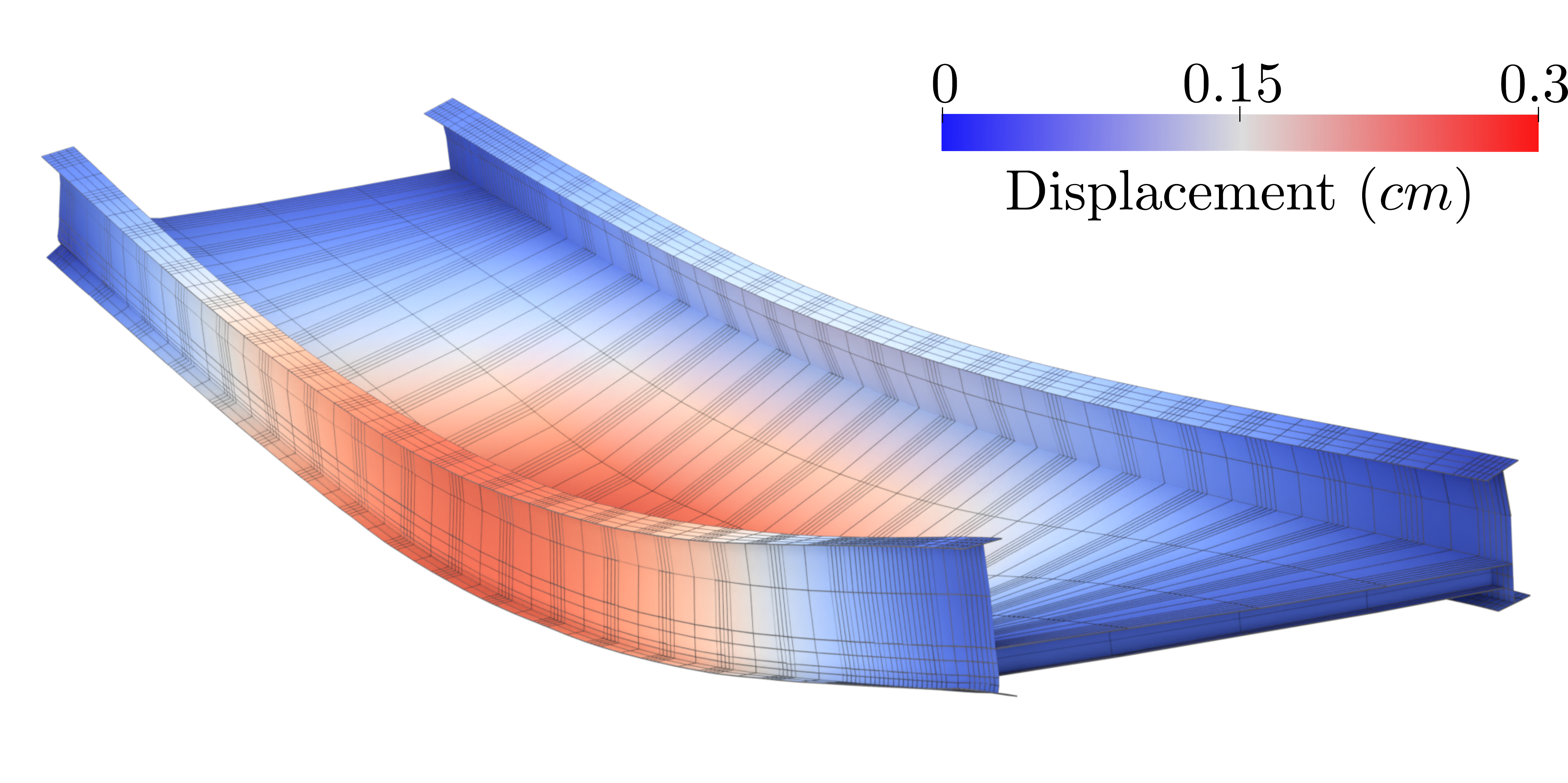}
	} 
	\subfloat[$t = 3 \, s$]{
		\includegraphics[width=0.3\textwidth]{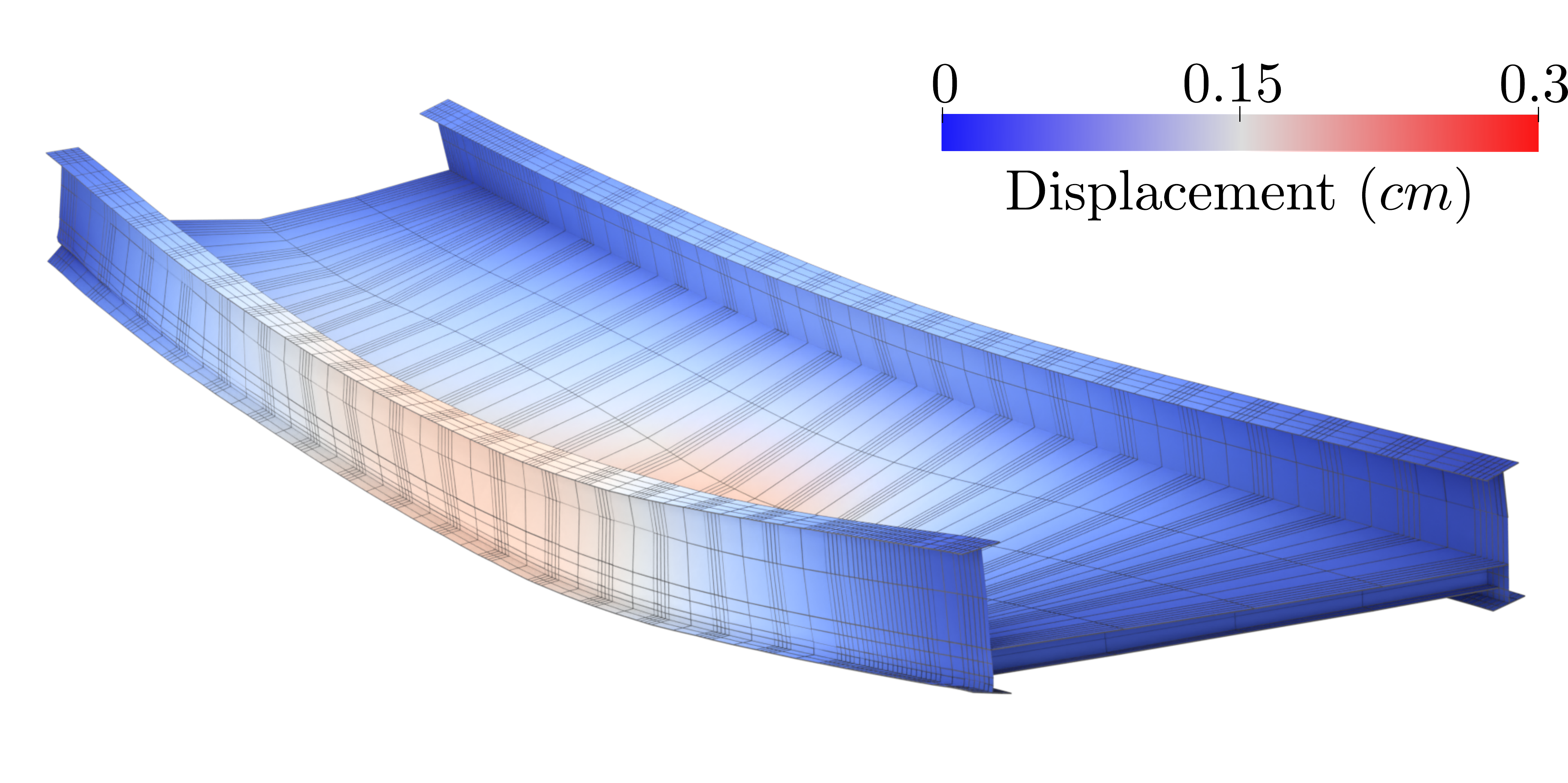}
	} \\
	\subfloat[$t = 1 .65\, s$ (at maximum deflection)]{
	\includegraphics[width=0.3\textwidth]{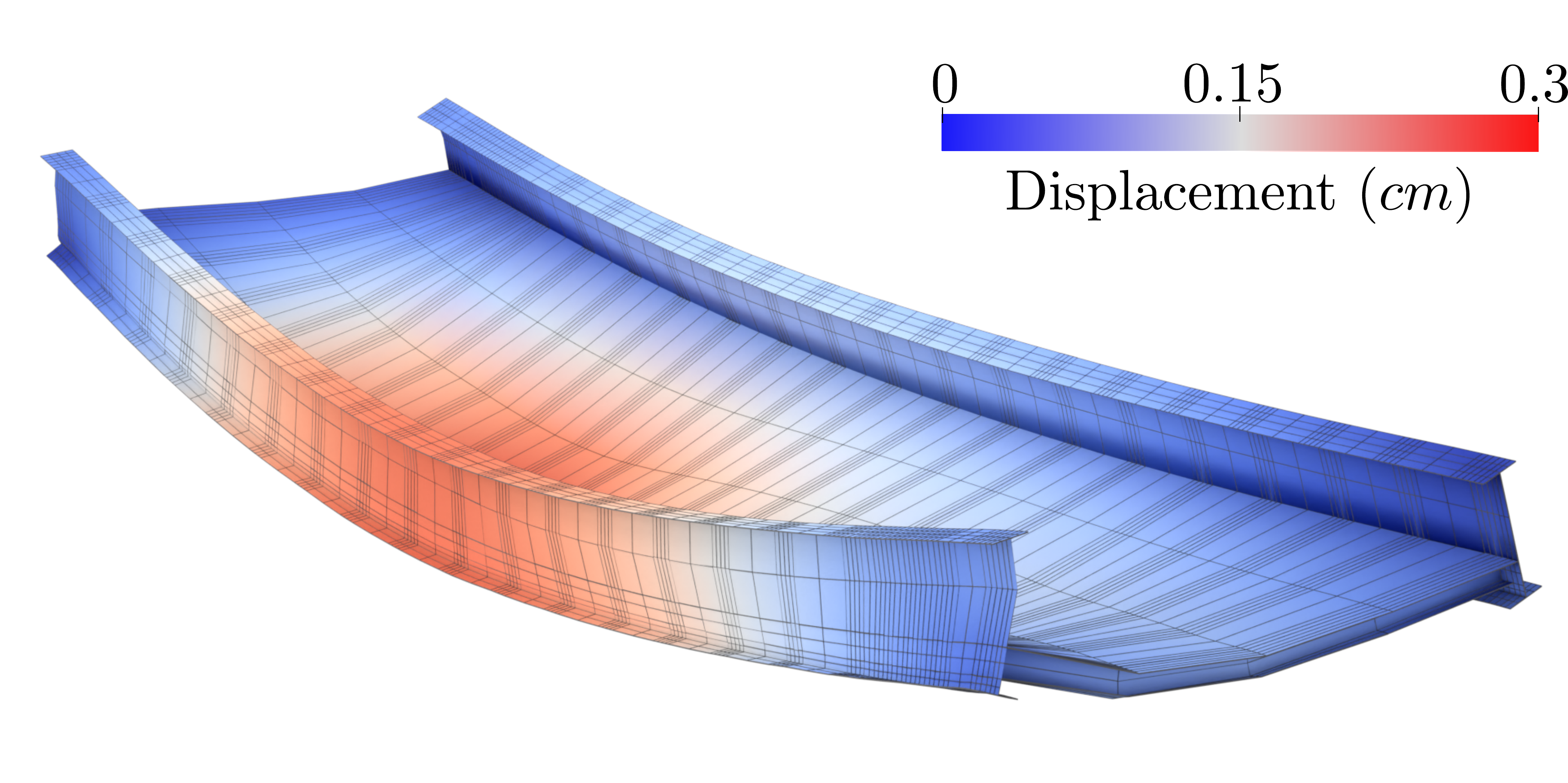} \label{fig:maxDeflection}
	}
	\subfloat[$t = 1 .75\, s$ (at maximum load)]{
	\includegraphics[width=0.3\textwidth]{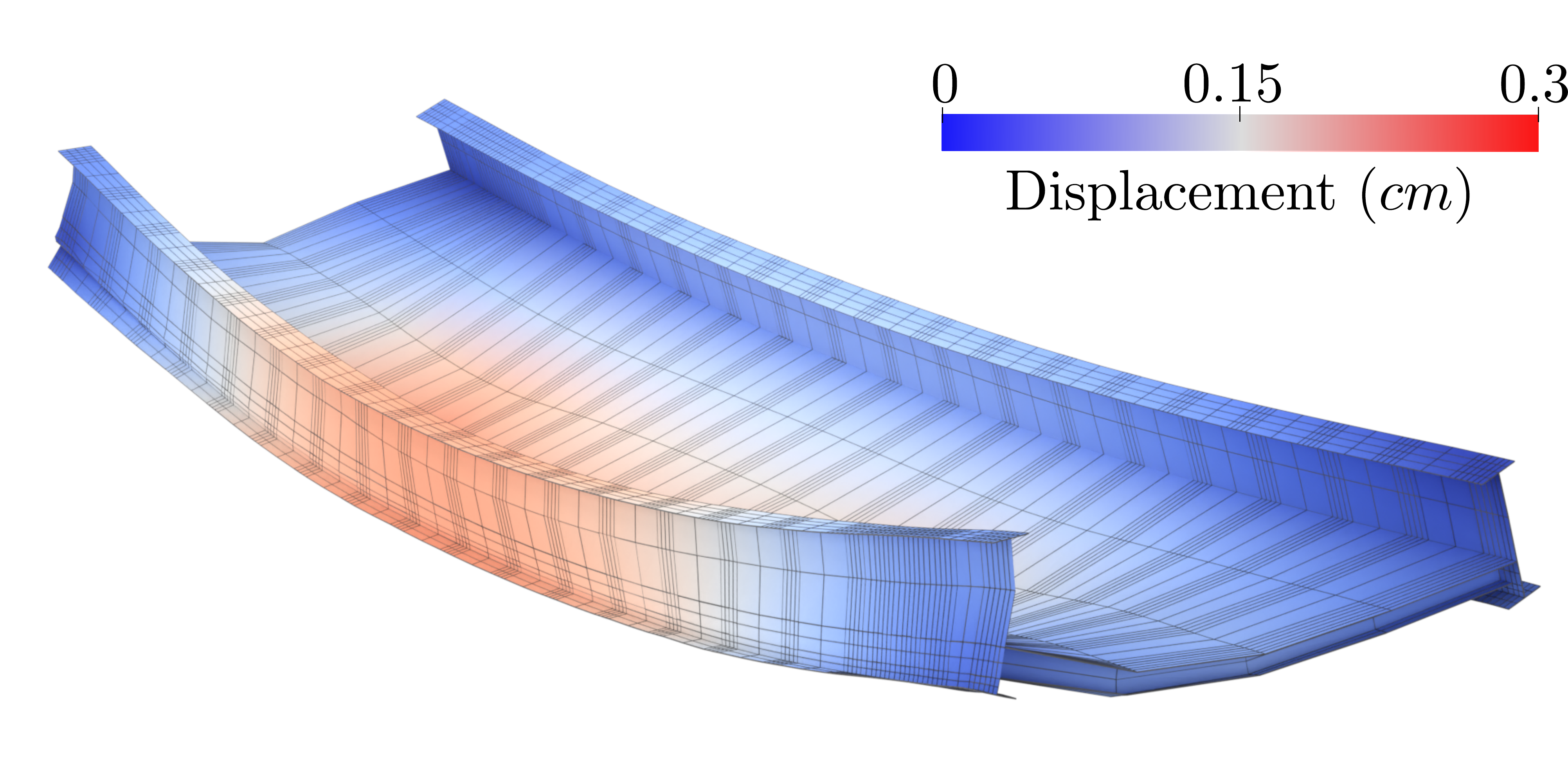} \label{fig:maxLoadDeflection}
	}
	\caption{Total train loading applied to the superstructure and its deflection at five distinct time instances} \label{fig:strainMesh4} 
\end{figure}

%
Throughout this section, a train of type T1 with four carriages heading from north to south on the east track is considered. Specifically, the analysis and discussion focuses on the axial strains along the top and bottom flanges of the east main I-beam. In the FE model the axial strains are obtained by multiplying the nodal displacement vector~$\ary u$ with the projection matrix $\ary P$ which discretises the symmetric gradient operator.  
The total loading depends on the number of train wheels on the bridge and is depicted in Figure~\ref{fig:forceOverTime}. The deflected superstructure at five distinct time instances is shown in Figure~\ref{fig:strainMesh4}. 
%

\subsection{Mesh convergence and probabilistic FE analysis \label{sec:forward}}
\begin{figure} []
	\setlength{\fboxsep}{0pt}%
	\setlength{\fboxrule}{0pt}%
	\centering
	\subfloat[FE strains $\ary P \ary u$ for meshes M1 and M2]{
		\includegraphics[width=0.46\textwidth]{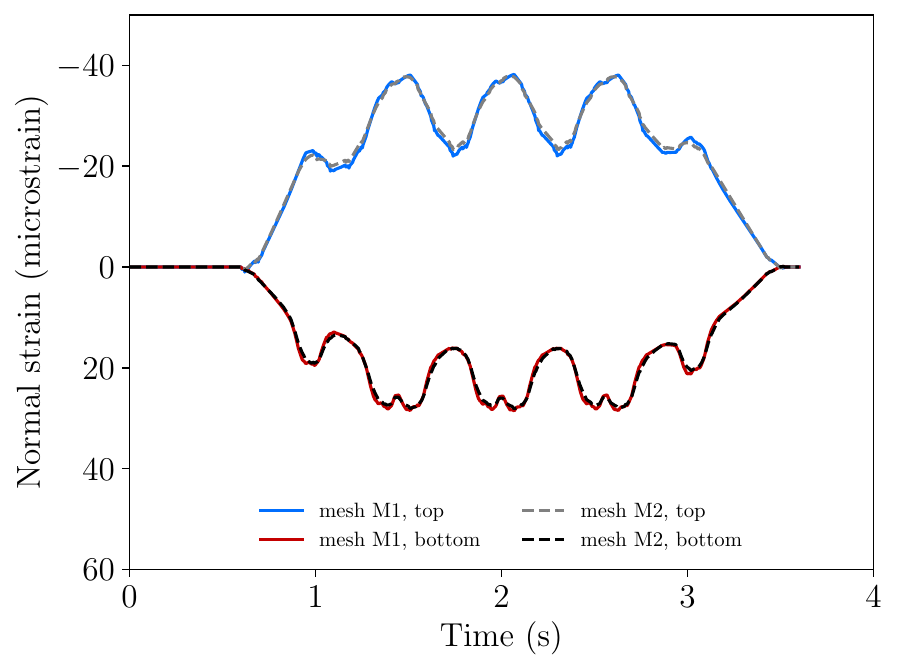} \label{fig:strainMeshes}
	} \quad
	\subfloat[Measured strains $\ary y$ and FE strains $\ary P \ary u$ for mesh M1]{
		\includegraphics[width=0.46\textwidth]{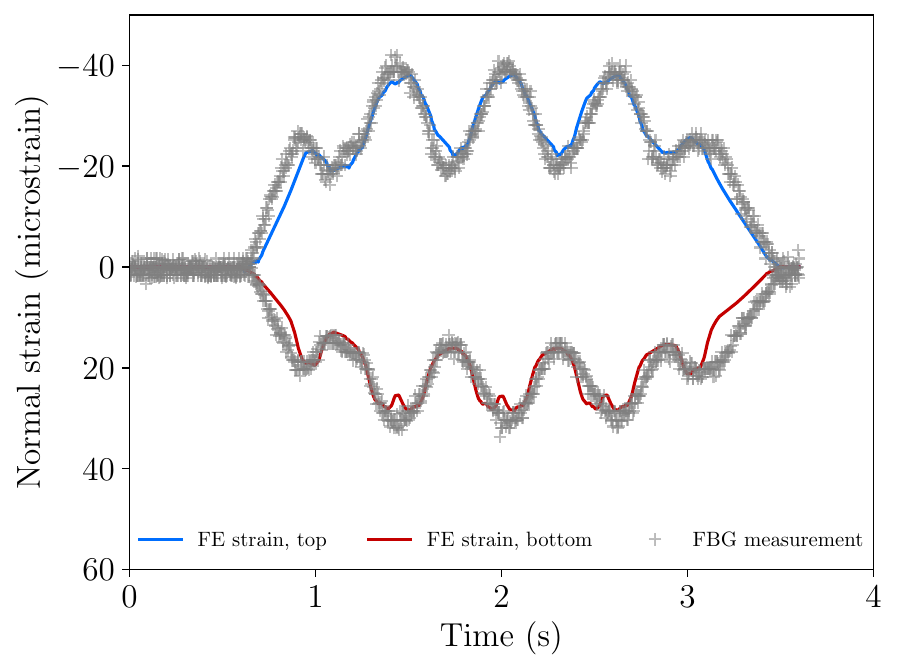} \label{fig:strainOverTime}
	} 
	\caption{Normal strains along the top and bottom flanges of the east main I-beam at the midspan\label{fig:forceStrain}} 
\end{figure}

%
The axial FE strains at the midspan of the east I-beam obtained using meshes M1 and  M2 are compared in Figure~\ref{fig:strainMeshes}. The results for both meshes are in close agreement; therefore the coarser mesh M1 is chosen for subsequent computations. The large-scale oscillations in the strains result from the changing number of train wheels on the bridge. The small-scale oscillations in the lower flange are due to the relative position of the wheel with respect to the relatively stiff transverse I-beams and the less stiff reinforced concrete deck.

The measured and computed FE strains at the midspan of the main I-beams are compared in Figure~\ref{fig:strainOverTime}. Each measured strain is depicted with a `{\small +}'. Evidently there is a significant scatter in the measured strains between neighbouring measurement points. Furthermore, the  offset between the measured and FE strains in time indicates a mismatch between the actual and the assumed train speed. The actual train speed and position corresponding to the measured strains are unknown. The time coordinate for the FE strains is chosen so that the offset between the two sets of strains is minimised.  Although the FE strains are smaller than the measured strains, their overall oscillation pattern is in close agreement. In comparison to the measured strains, the bottom flange FE strains have multiple small dips which is again related to the modelling of the connection between the main and transverse I-beams and the concrete deck. 

The random loading~$\vec r(\vec x)$ representing the uncertainties in the train loading and its application to the superstructure yields according to~\eqref{eq:randomU} the prior FE strain density $p(\ary P \ary u) = \set N (\ary P \overline{ \ary u} , \, \ary P \ary C_{\ary u} \ary P^\trans) $.  The mean and the 95\% confidence region of the strains computed at 20 sensor locations along the top and bottom flanges of the east main I-beam (i.e. 40 sensor locations in total) are shown in Figure~\ref{fig:pUgivY501}. The plotted curves are non-smooth because the strains are linearly interpolated between the sensor locations. The actual FE strains are smooth and quadratic within each element~\citep{zhang2018subdivision}. 

\subsection{Statistical finite element analysis \label{sec:statFEMresults}}
%

%
\subsubsection{Inference of true strain response using all sensors \label{sec:trueSystemFull}}
%

%
\begin{table}[t]
	\centering
	\begin{tabular} { l c c  }
		\toprule  
		parameter & value   \\
		\midrule
		$\sigma_r$ & 1000 Pa \\
		$\ell_r$ & 1 m \\
		$\sigma_e$ & 1 microstrain \\
		MCMC samples & 20000 \\
		\bottomrule
	\end{tabular}
	\caption[]{Covariance and algorithmic parameters used in Section~\ref{sec:statFEMresults} \label{tab:params}}
\end{table}

The measured strains are used to predict the `true' axial strains along the east main I-beam. The posterior FE strain density~\mbox{$p(\ary P \ary u | \ary y)$} and the true strain density~\mbox{$p( \ary z | \ary y)$} conditioned on the measured strains~$\ary y$ are determined by evaluating~\eqref{eq:posteriorU} and~\eqref{eq:posteriorZ}. Evidently, the true strain response~$\ary z$ is unknown and only the measured strain response including some measurement error, i.e. \mbox{$\ary y = \ary z + \ary e$}, at the $n_y$ locations is known. At each time instant~$t_k$,  the posterior FE strain density~\mbox{$p(\ary P \ary u_k | \ary y_k)$} and the true strain density~\mbox{$p( \ary z_k | \ary y_k)$} are computed from the measured strains $\ary y_{k}$. For simplicity, the subscript~$k$ is omitted here and in the following.

\begin{figure*} []
	\setlength{\fboxsep}{0pt}%
	\setlength{\fboxrule}{0pt}%
	\centering
	\subfloat[$p( \ary Y | \rho )$]{
		\includegraphics[width=0.31\textwidth]{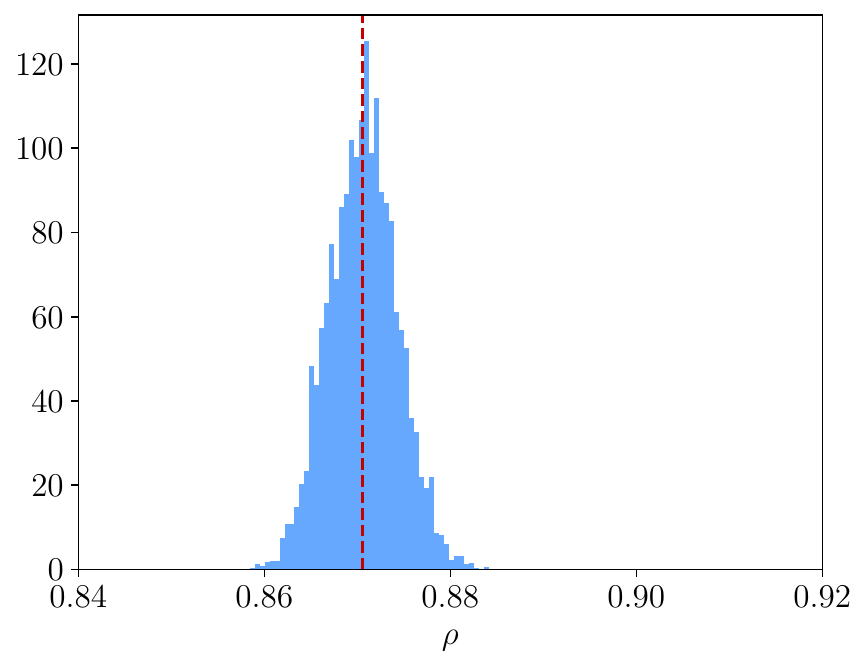} 
	} 
	\subfloat[$p( \ary Y | \sigma_d )$]{
		\includegraphics[width=0.31\textwidth]{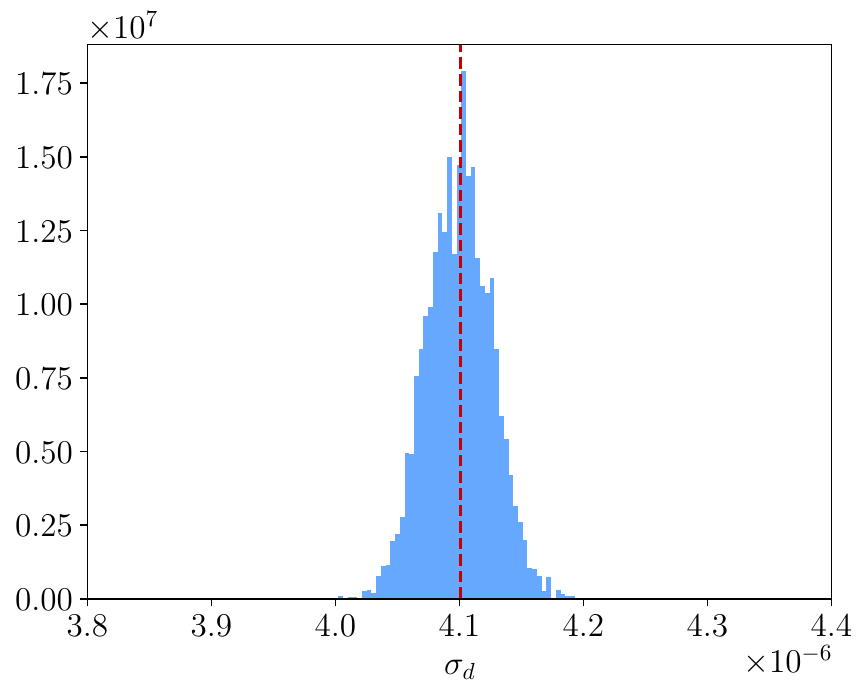}
	} 
	\subfloat[$p( \ary Y | \ell_d )$]{
		\includegraphics[width=0.31\textwidth]{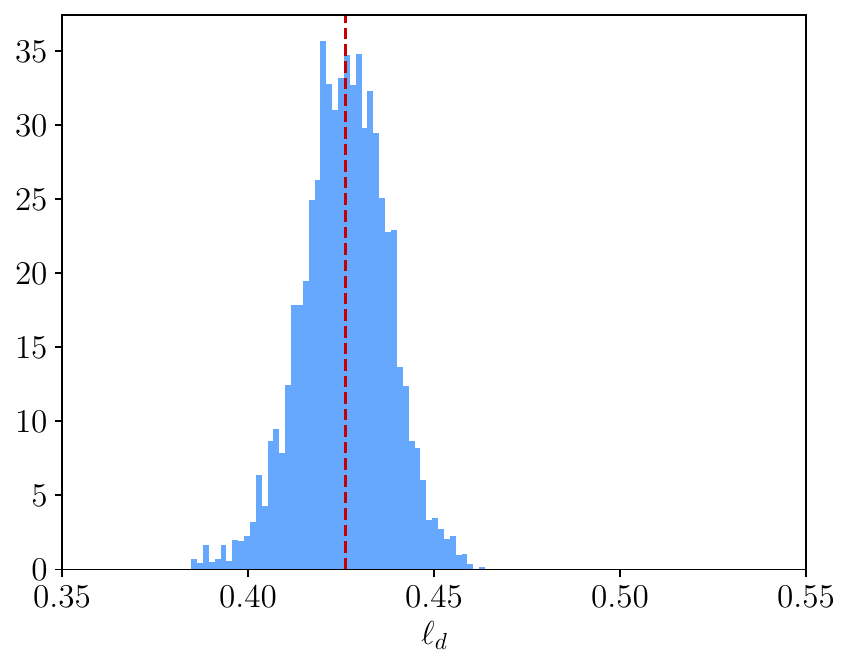}
	} 
	\caption{Normalised histogram of the marginal likelihood~$p(\ary Y | \ary w)$ for $n_y=40$ and~$n_o = 501$ sampled with MCMC. The red dashed lines indicate the empirical mean~$\ary w^*$ of the samples} \label{fig:hParam501} 
\end{figure*}

Before computing the posterior densities, the marginal likelihood~\mbox{$p(\ary Y | \ary w)$} given by~\eqref{eq:marginalY} is sampled to obtain the point estimate~$\ary w^*$. Recall that the hyperparameter vector~$\ary w$ consists of the scaling parameter~$\rho$ and the covariance kernel parameters~$\sigma_d$ and~$\ell_d$ of the mismatch error~$\ary d$. The strain recordings from all $n_y = 40$ FBG sensors of the east main I-beam are included in the observation (measurement) matrix \mbox{$\ary Y \in \mathbb R^{n_y \times n_o}$}. The sensor locations are as specified in Figure \ref{fig:bridgeschematic}. Only data recorded when a train is present on the bridge (i.e., \mbox{$1~\text{s} \leq t \leq 3~\text{s}$} yielding $n_o = 501$ readings for each sensor) are included. The standard deviation $\sigma_e =  1~\text{microstrain}$ of the measurement error~$\ary e$ is estimated from the the strain recordings within \mbox{$0~\text{s} \leq t \leq 0.5~\text{s}$} prior to the arrival of the train (see Figure \ref{fig:strainOverTime}). The marginal likelihood~$p(\ary Y | \ary w)$  is sampled with a standard MCMC algorithm.  In Figure~\ref{fig:hParam501} the normalised histograms for $p(\ary Y | \rho )$, $p(\ary Y | \sigma_d)$, and $p(\ary Y| \ell_d )$ obtained with 20000 iterations and an acceptance ratio of 0.283 are depicted. Convergence is ensured by checking the trace plot of the samples and monitoring the empirical mean and standard deviation of the samples. As mentioned before, the empirical mean of the MCMC samples is used as a point estimate~$\ary w^*$, see Table~\ref{tab:hParamStat501}.  The unimodal nature of the histograms in Figure~\ref{fig:hParam501} and the relatively small standard deviations indicate that the available data is sufficient to identify the hyperparameters with a high certainty. Lack of sufficient data usually leads to multimodal histograms and large standard deviations.
%


%
\begin{table}[t]
	\centering
	\begin{tabular} { l c c  }
		\toprule  
		hyperparameter & value   \\
		\midrule
		$\rho$ & 
		$0.8706 \pm 0.0037 $
		\\
		$\sigma_d$ (microstrain)& 
		$ 4.0998 \pm 0.0263 $ 
		\\
		$\ell_d$ (meter)& 
		$ 0.4261 \pm 0.0116$ 
		\\
		\bottomrule
	\end{tabular}
	\caption[]{Empirical mean and standard deviation of the hyperparameters. \label{tab:hParamStat501}}
\end{table}

\begin{figure*} []
	\setlength{\fboxsep}{0pt}%
	\setlength{\fboxrule}{0pt}%
	\centering
	\subfloat[$t = 1$s]{
		\includegraphics[width=0.32\textwidth]{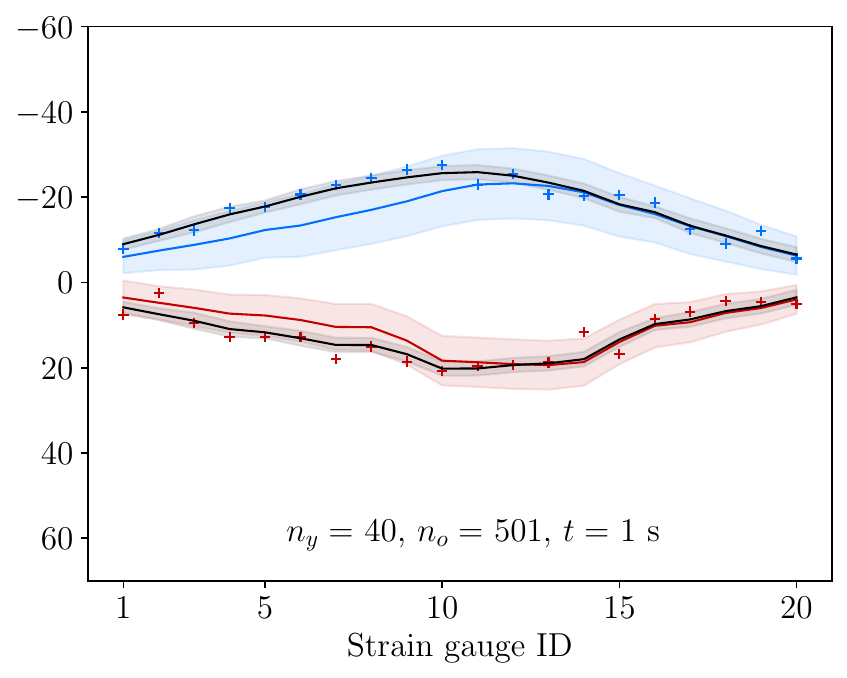} 
	} 
	\subfloat[$t = 2$s]{
		\includegraphics[width=0.32\textwidth]{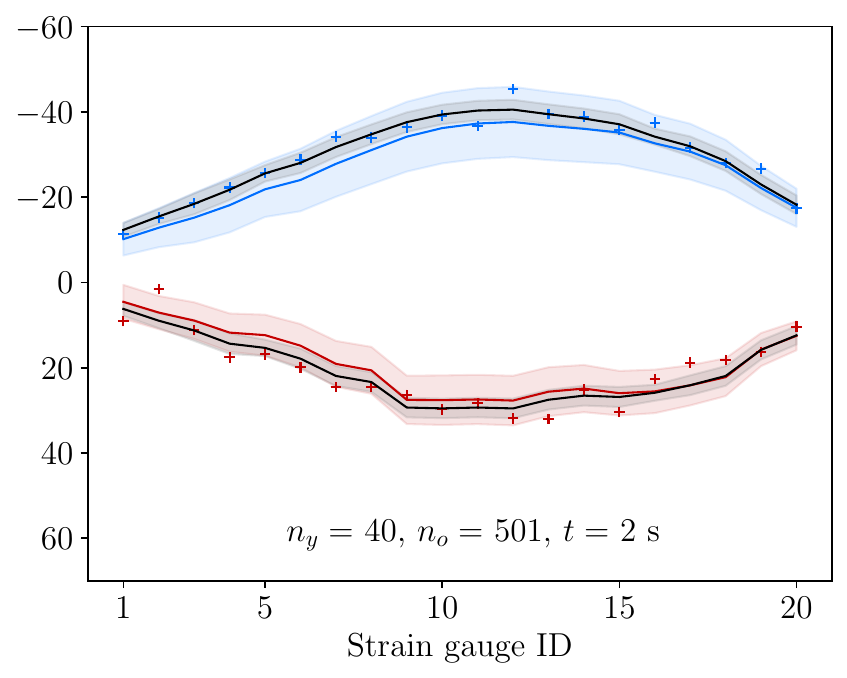}
	} 
	\subfloat[$t = 3$s]{
		\includegraphics[width=0.32\textwidth]{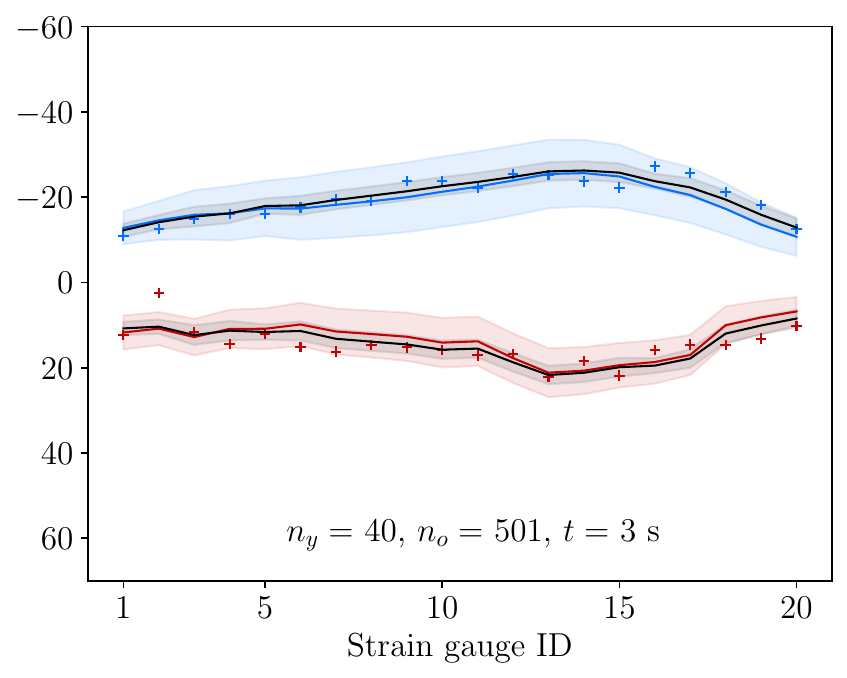}
	} \\

	\caption{Posterior FE strains~$p(\ary u |\ary y)$  conditioned on the measured strains (+) along the east main I-beam. The blue and red lines represent the mean~$\ary P \overline{ \ary u }$ of the prior along the top and bottom flanges, respectively, and the black lines the conditioned mean~$\rho \ary P \overline {\ary u}_{|\ary y}$. The shaded areas denote the corresponding $95\%$ confidence regions. The unit of the vertical axis is microstrain} \label{fig:pUgivY501} 
\end{figure*}

\begin{figure*} []
	\setlength{\fboxsep}{0pt}%
	\setlength{\fboxrule}{0pt}%
	\centering
	\subfloat[$t = 1$s]{
		\includegraphics[width=0.32\textwidth]{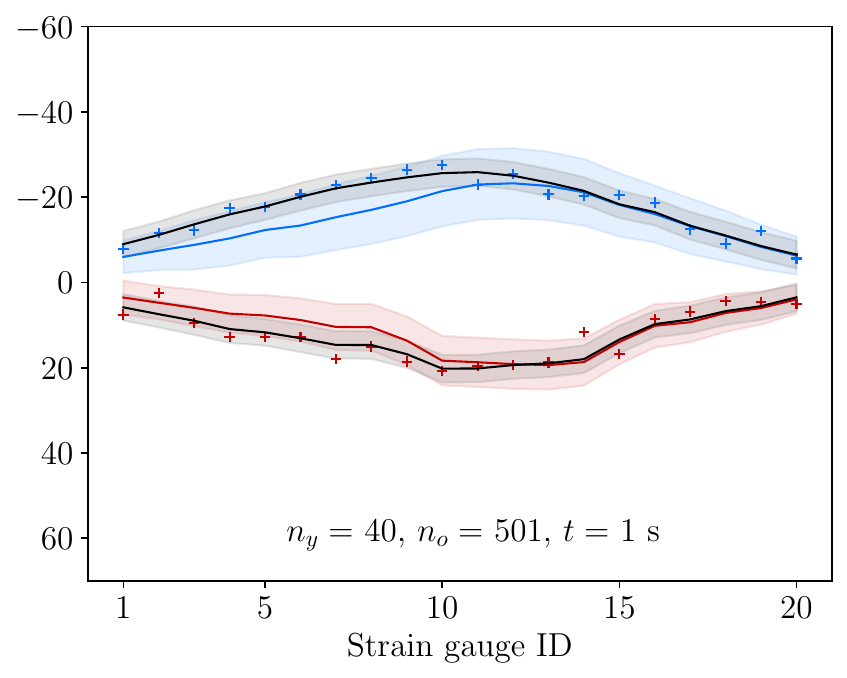} 
	} 
	\subfloat[$t = 2$s]{
		\includegraphics[width=0.32\textwidth]{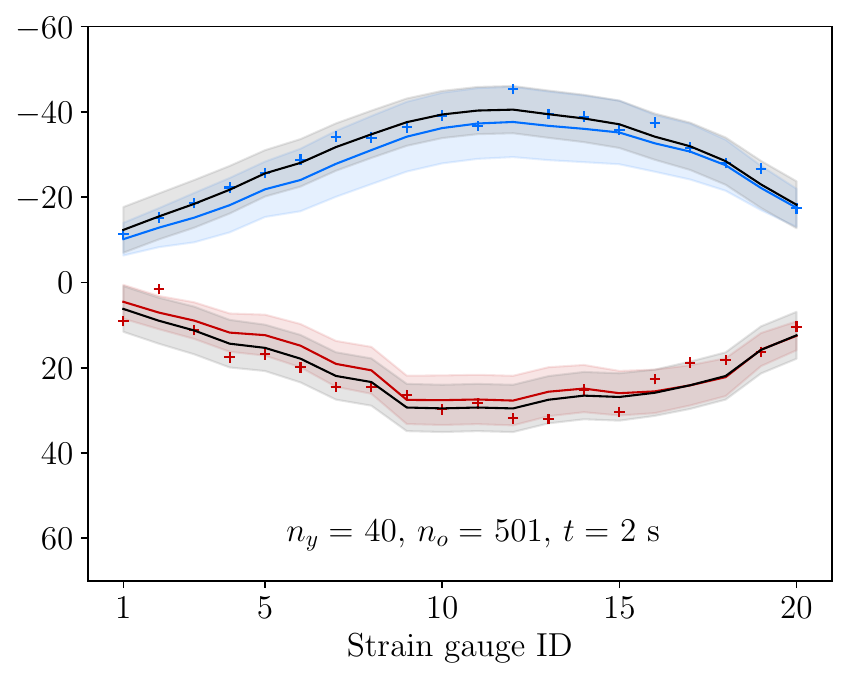}
	} 
	\subfloat[$t = 3$s]{
		\includegraphics[width=0.32\textwidth]{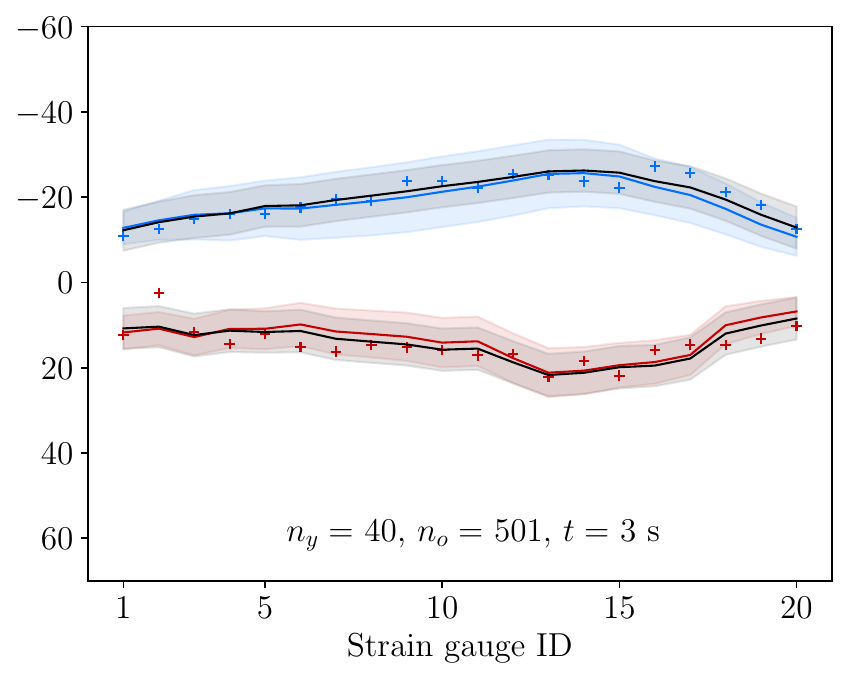}
	} 
	\caption{Inferred true strain density~$p(\ary z | \ary y)$  conditioned on strains (+) measured along the east main I-beam. The blue and red lines represent the mean~$\ary P \overline{\ary u }$ of the prior along the top and bottom flanges, respectively, and the black lines the conditioned mean~$\overline {\ary z}_{|\ary y}$. The shaded areas denote the corresponding $95\%$ confidence regions. The unit of the vertical axis is microstrain} \label{fig:pZgivY501} 
\end{figure*}

In Figure~\ref{fig:pUgivY501} the prior and posterior densities of the axial FE strains and the measured strains at \mbox{$t = 1~\text{s}$}, \mbox{$t=2~\text{s}$}, and \mbox{$t=3~\text{s}$} are plotted. For the FE strains, in addition to the mean, the corresponding 95\% confidence regions are plotted. Observe how the mean of the posterior FE strain lies much closer than the prior FE strain to the measured strains. It is also evident that the uncertainty in the prior FE strains is significantly reduced by conditioning them on the sensor measurements. In the corresponding Figure~\ref{fig:pZgivY501} the posterior densities of the true system response at \mbox{$t = 1~\text{s}$},  \mbox{$t=2~\text{s}$}, and \mbox{$t=3~\text{s}$} are shown. The respective confidence regions encompass all the strain measurements and their mean is always visually very close to the mean of the data. These results clearly demonstrate that the FE strains conditioned on the measured strains (i.e. black curves in Figure~\ref{fig:pUgivY501}) provide an improvement of strain prediction over the unconditioned FE strain results. Moreover, the 95\% confidence regions provide a quantifiable method for identifying anomalous sensor readings. As a further note, observe that the confidence interval of $p(\ary z | \ary y)$ in Figure  \ref{fig:pZgivY501} is wider than the confidence interval of  $p(\ary u | \ary y)$ in Figure \ref{fig:pUgivY501} due to the contribution of the inferred mismatch term.


\subsubsection{Inference of true strain response using a reduced number of sensors  \label{sec:trueSystemRed}}
%
Given the considerable efforts and costs associated with instrumenting operational structures, methods for optimising both number and location of sensing points are critical.  The statFEM approach allows for both data and physics-informed prediction even in situations where very limited measurement data are available.   In the following, a reduced number of measurements are considered to obtain the posterior FE strain density~$p(\ary P \ary u | \ary y)$ and the true strain density~$p(\ary z | \ary y)$. A reduced number of sensors $n_y$ and readings per sensor $n_o$ are considered. The purpose of this study is to confirm empirically the convergence of the two posteriors with increasing number of data and to determine the minimum number of  data required for an acceptable estimate. The use of fewer data has advantages in terms of reduced instrumentation, increased numerical efficiency and the corresponding reduction in effort spent on data analysis and interpretation.

Three numbers of sensing points \mbox{$n_y \in \{ 40, \, 20, \, 10 \} $} along the east main I-beam are considered. In each case half of the sensing points are along the top and the other half along the bottom flange. Recall that the total number of sensors installed along the top and bottom flanges of each main I-beam are 40 (20 top and 20 bottom).  For each~$n_y$ a subset of FBG sensors with the IDs given in Table~\ref{tab:gaugeID} is selected. 
\begin{table}[]
	\centering
	\begin{tabular} { l c c c }
		\toprule  
		& gauge ID   & top/bottom flange\\
		\midrule
		$n_y = 40$ & $1, 2, 3, \ldots , 20$ & top \& bottom
		\\
		$n_y = 20$& $ 1, 3, 5, \ldots , 19 $ &  top \& bottom
		\\
		$n_y = 10$& $ 3, 7, 11, 15, 19$  &  top \& bottom
		\\
		\bottomrule
	\end{tabular}
	\caption[]{Strain gauge ID for computations with $n_y = \{40, 20, 10\}$ sensors located on the east main I-beam. \label{tab:gaugeID}}
\end{table}
Furthermore, the number of  strain measurements~$n_o$ is altered by selecting four different time intervals \mbox{$\Delta t = \{1/250~\text{s}, \, 1/50~\text{s}, \,  1/25~\text{s}, \, 1~\text{s}\}$} between the measurements within the observation window \mbox{$1~\text{s} \leq t \leq 3~\text{s}$}. The respective number of measurements  are $n_o = \{ 501, \, 101, \,  51, \, 3\}$. As in the preceding section, the marginal likelihood~$p(\ary Y | \ary w)$ is sampled using MCMC. For each combination of $n_y$ and $n_o$, a total of 20000 MCMC samples are generated with an average acceptance ratio of 0.326. 

\begin{figure*} []
	\setlength{\fboxsep}{0pt}%
	\setlength{\fboxrule}{0pt}%
	\centering
	\subfloat[$p(\ary Y | \rho  ), \, n_y = 20$]{
		\includegraphics[width=0.32\textwidth]{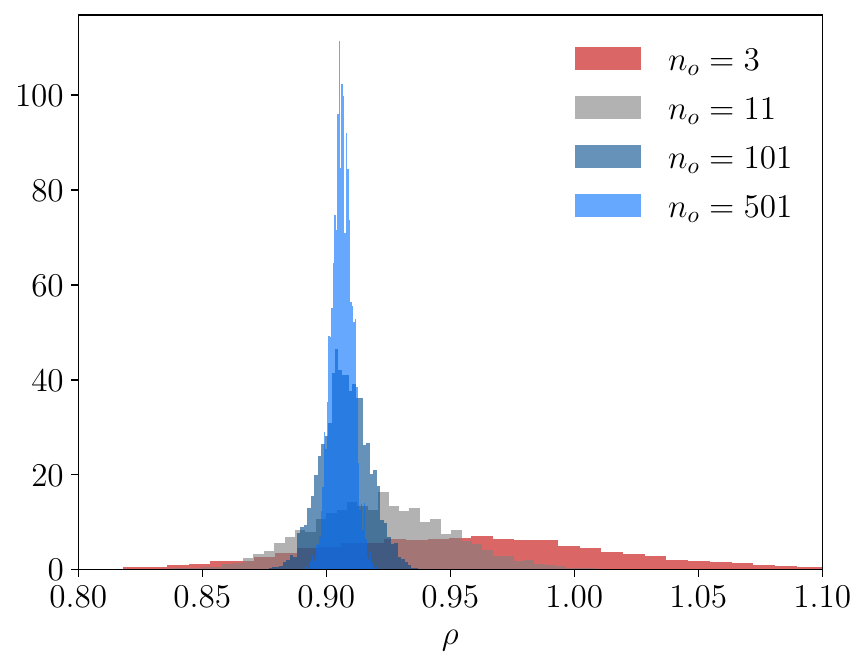} 
	} 
	\subfloat[$p( \ary Y | \sigma_d ), \, n_y = 20$]{
		\includegraphics[width=0.32\textwidth]{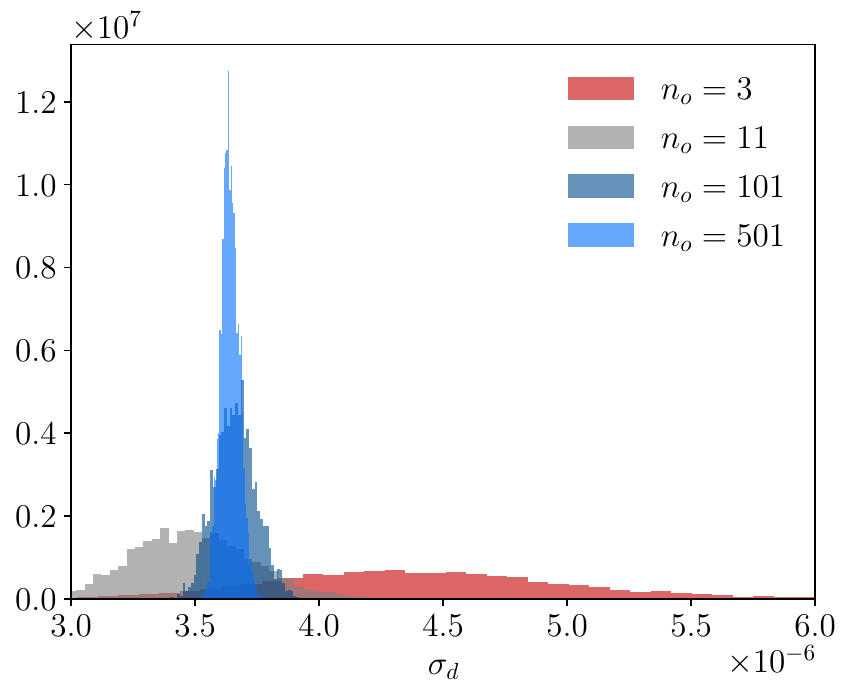}
	} 
	\subfloat[$p(\ary Y | \ell_d ), \, n_y = 20$]{
		\includegraphics[width=0.32\textwidth]{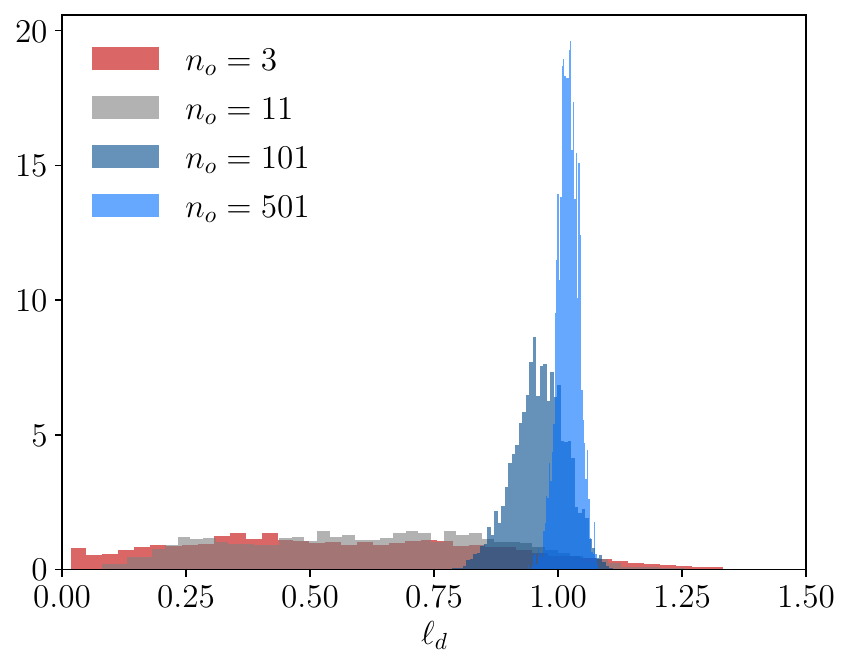}
	} 
	\caption{Normalised histogram of the marginal likelihood~$p(\ary Y | \ary w)$ for~$n_y=20$ and~$n_o \in \{3, \, 11,\, 101, \, 501 \}$ sampled with MCMC} \label{fig:hParam} 
\end{figure*}

The normalised histograms for \mbox{$p(\ary Y | \rho )$}, \mbox{$p(\ary Y | \sigma_d )$}, and \mbox{$p(\ary Y | \ell_d )$} for \mbox{$n_y = 20$} are shown in Figure~\ref{fig:hParam}. 
It is apparent that the standard deviations become significantly smaller with increasing number of readings~$n_o$. The empirical mean and standard deviation for all combinations of $n_y$ and $n_o$ are given in the Tables \ref{tab:rho}, \ref{tab:sigmaD}, and \ref{tab:ellD}. In almost all cases the standard deviation of the samples becomes smaller with increasing number of readings. Moreover, the standard deviation becomes smaller when data from more sensors are considered (larger $n_y$). For the scaling parameter $\rho$, the difference between the means obtained with \mbox{$n_o = 501$} and \mbox{$n_o=3$} is around $5\%$ and decreases as additional readings are incorporated. For the mismatch parameter $\sigma_d$, the difference can be as high as $25 \%$, and then decreases to 2--8\% depending on the number of sensors $n_y$. It is remarkable that with 50 times fewer readings, the difference of the estimated $\rho$ and $\sigma_d$ is still acceptable. The length scale parameter $\ell_d$, however, shows a larger difference between the two smaller number of sensing points $n_o =\{3,  \,11\}$ and the larger number of sensing points $n_o = 501$. This result, however, does not significantly affect the inferred posterior true strains.

\begin{table*}[]
	\centering
	\begin{tabular} { l c c c c }
		\toprule  
		& $n_y = 10$ & $n_y = 20$ & $n_y = 40$  \\
		\midrule
		$n_o = 3$ & 
		$1.0152 \pm 0.0717$ 
		& $ 0.9581 \pm 0.0593$ 
		& $ 0.9253 \pm 0.0548$ 
		\\
		$n_o = 11$ & 
		$ 0.9614 \pm 0.0325$ 
		& $ 0.9215 \pm 0.0281$ 
		& $ 0.8869 \pm 0.0251$ 
		\\
		$n_o = 101$ & 
		$ 0.9452 \pm 0.0103$ 
		& $ 0.9075 \pm 0.0093$ 
		& $ 0.8717 \pm 0.0080 $ 
		\\
		$n_o = 501$ & 
		$ 0.9457 \pm 0.0046$ 
		& $ 0.9062 \pm 0.0041$ 
		& $ 0.8706\pm 0.0036 $ 
		\\
		\bottomrule
	\end{tabular}
	\caption[]{Empirical mean and standard deviation of the hyperparameter~$\rho$. \label{tab:rho}}
\end{table*}
\begin{table*}[]	
	\centering
	\begin{tabular} { l c c c c }
		\toprule  
		& $n_y = 10$ & $n_y = 20$ & $n_y = 40$  \\
		\midrule
		$n_o = 3$ & 
		$4.9211 \pm 1.0585$ 
		& $ 4.4129 \pm 0.6042$ 
		& $ 4.6830 \pm  0.4241 $ 
		\\
		$n_o = 11$ & 
		$ 3.5436 \pm 0.3651$ 
		& $ 3.4951 \pm 0.2365$ 
		& $  3.9374 \pm 0.1742 $ 
		\\
		$n_o = 101$ & 
		$ 3.8462 \pm 0.1219 $ 
		& $ 3.6611\pm 0.0840$ 
		& $ 4.0894 \pm 0.0590$ 
		\\
		$n_o = 501$ & 
		$ 3.8524 \pm 0.0542$ 
		& $ 3.6405 \pm 0.0367$ 
		& $ 4.0998 \pm 0.0263$ 
		\\
		\bottomrule
	\end{tabular}
	\caption[]{Empirical mean and standard deviation of the hyperparameter~$\sigma_d$. \label{tab:sigmaD}}
\end{table*}
\begin{table*}[]
	\centering
	\begin{tabular} { l c c c c }
		\toprule  
		& $n_y = 10$ & $n_y = 20$ & $n_y = 40$  \\
		\midrule
		$n_o = 3$ & 
		$ 0.5657 \pm 0.4098$ 
		& $ 0.5634\pm 0.3060$ 
		& $ 0.3766 \pm 0.1442 $ 
		\\
		$n_o = 11$ & 
		$ 0.4663 \pm 0.3378$ 
		& $ 0.6099 \pm 0.2592$ 
		& $ 0.2508 \pm 0.1250$ 
		\\
		$n_o = 101$ & 
		$ 0.3736 \pm 0.1482$ 
		& $ 0.9628 \pm 0.5261$ 
		& $ 0.4335 \pm 0.0235 $ 
		\\
		$n_o = 501$ & 
		$ 0.4731\pm 0.1005$ 
		& $ 1.0201 \pm 0.0208$ 
		& $ 0.4260 \pm 0.0116$ 
		\\
		\bottomrule
	\end{tabular}
	\caption[]{Empirical mean and standard deviation of the hyperparameter~$\ell_d$. \label{tab:ellD}}
	\vspace{1.2em}
\end{table*}

\begin{figure*} []
	\setlength{\fboxsep}{0pt}%
	\setlength{\fboxrule}{0pt}%
	\centering
	\includegraphics[width=0.31\textwidth]{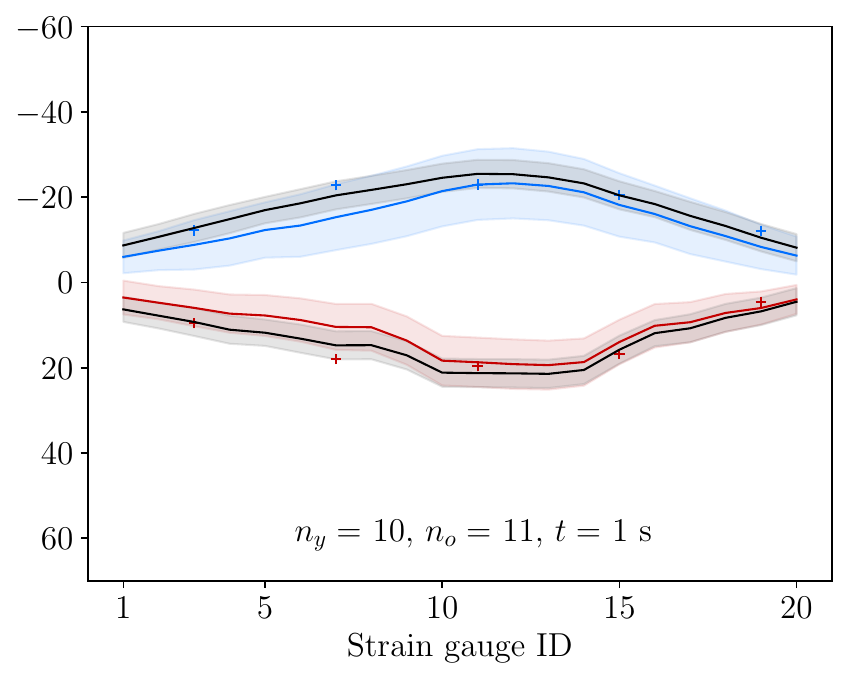} 
	\includegraphics[width=0.31\textwidth]{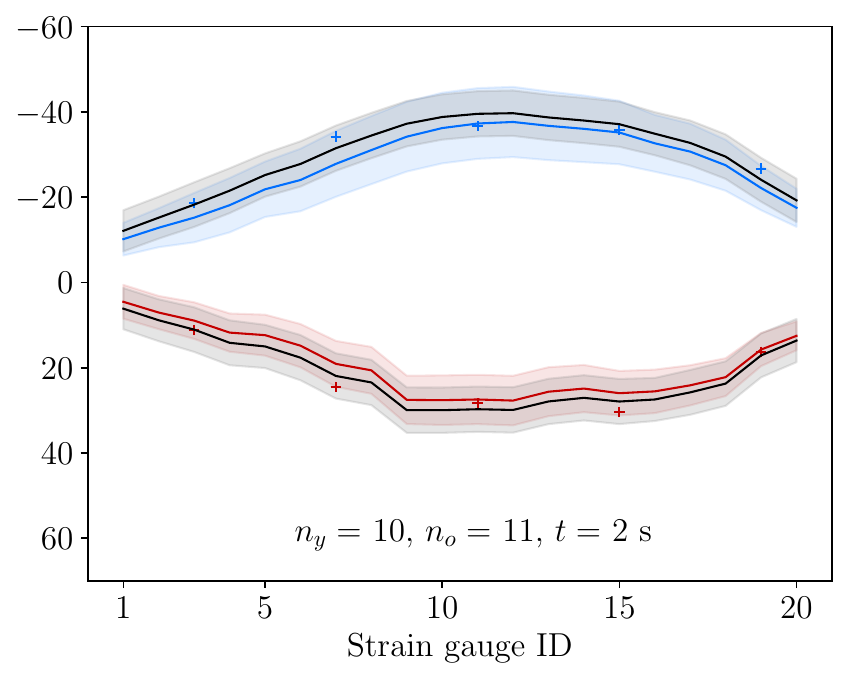}
	\includegraphics[width=0.31\textwidth]{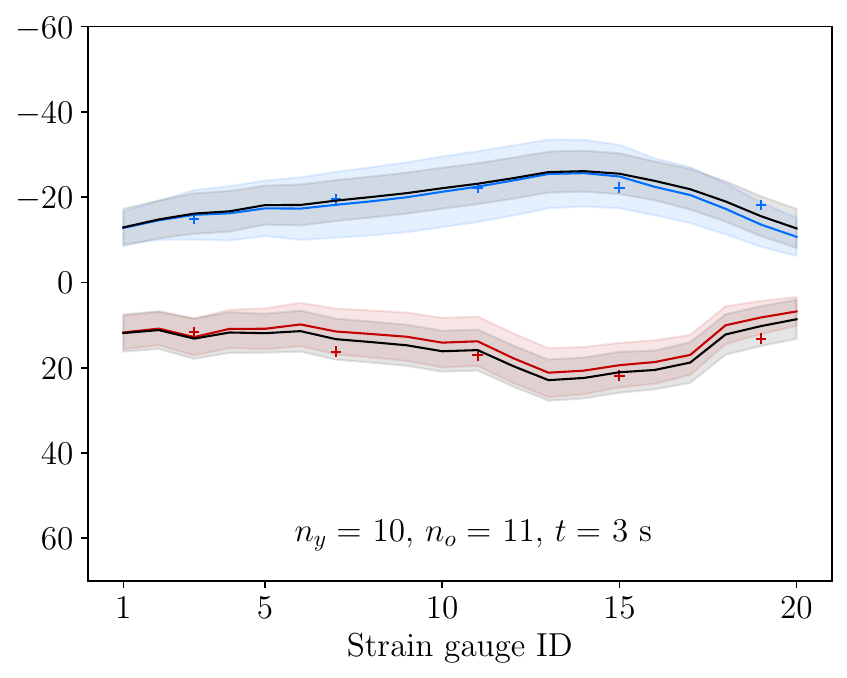}
	\\
	\includegraphics[width=0.31\textwidth]{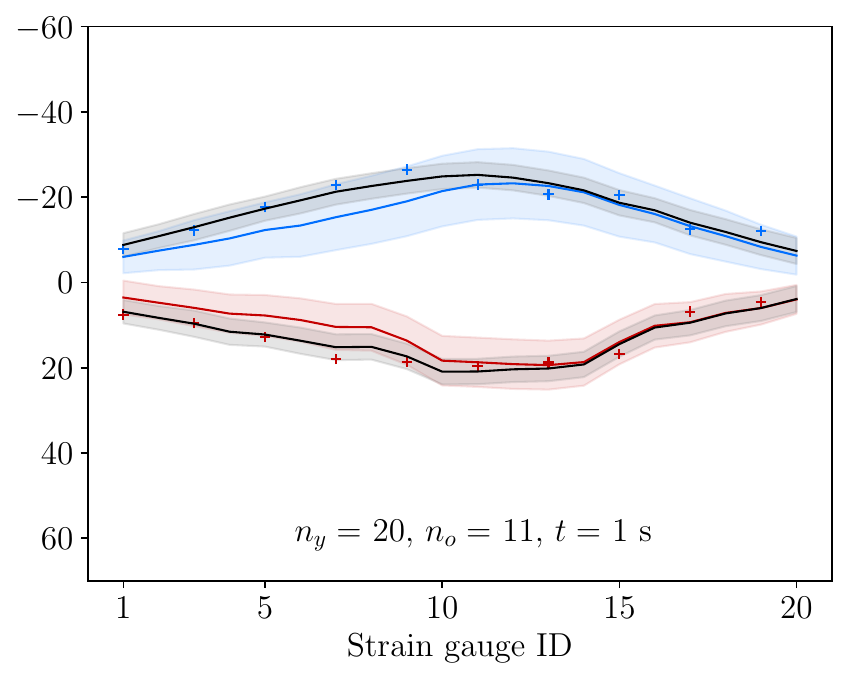} 
	\includegraphics[width=0.31\textwidth]{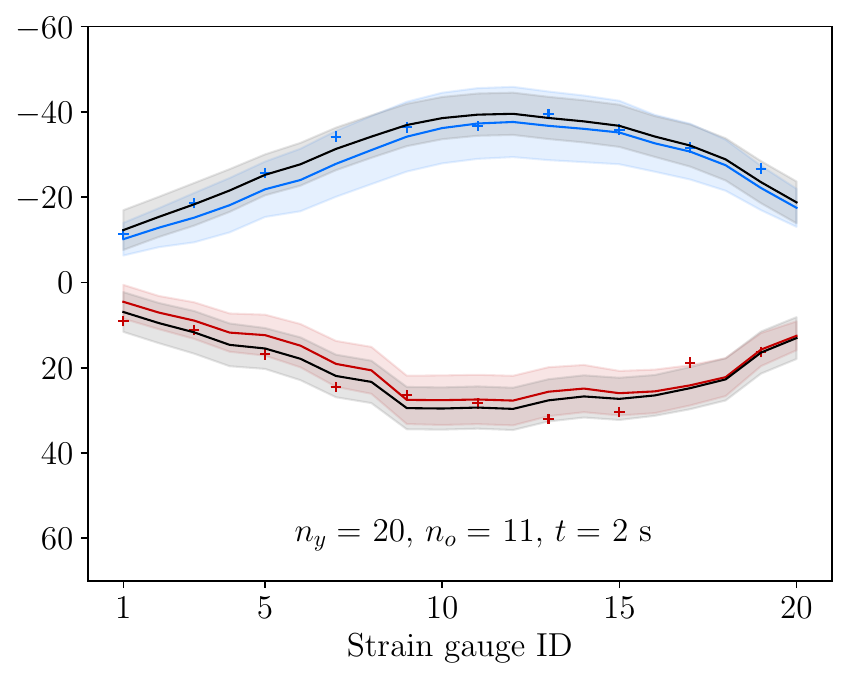}
	\includegraphics[width=0.31\textwidth]{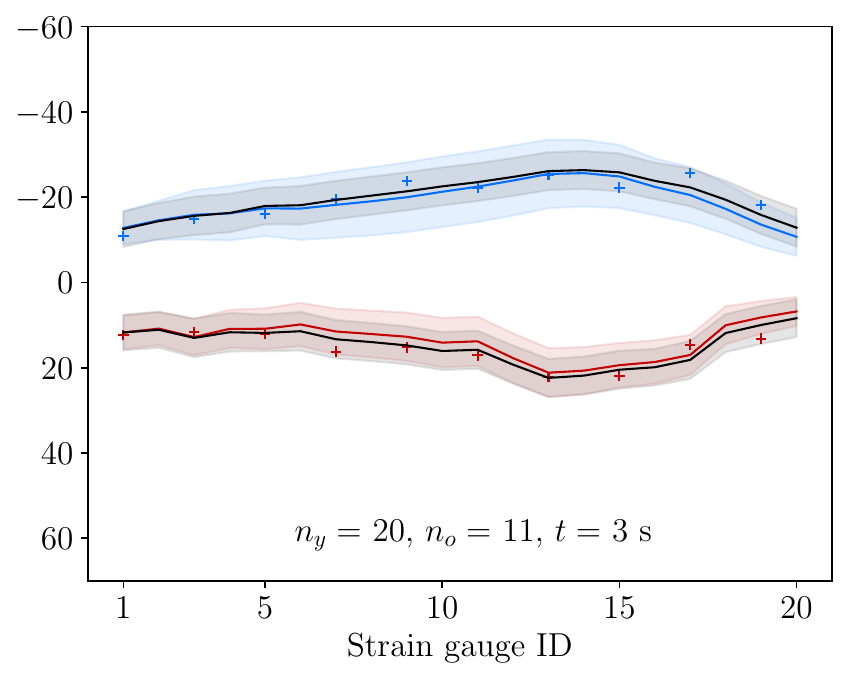}
	\\
	\includegraphics[width=0.31\textwidth]{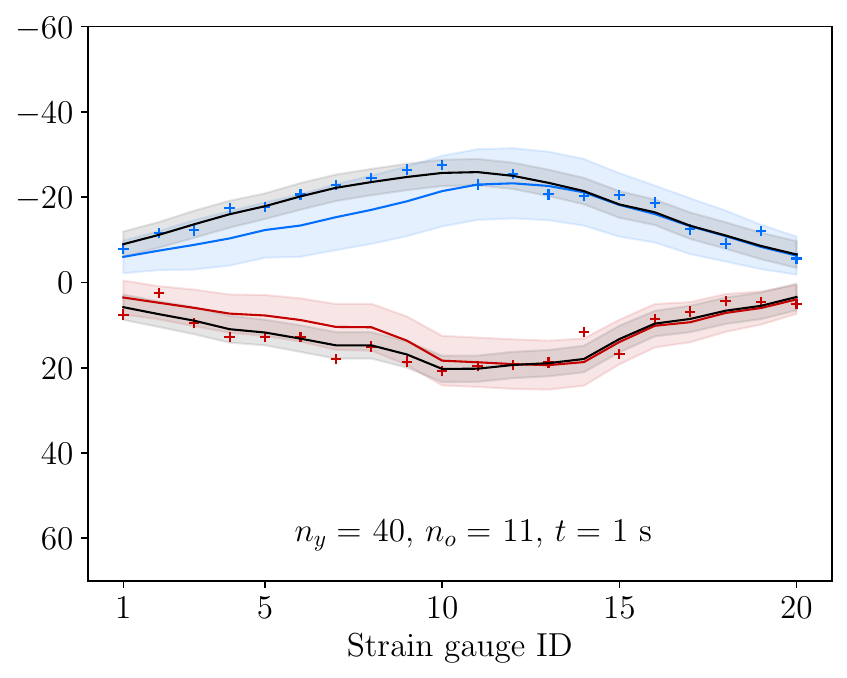} 
	\includegraphics[width=0.31\textwidth]{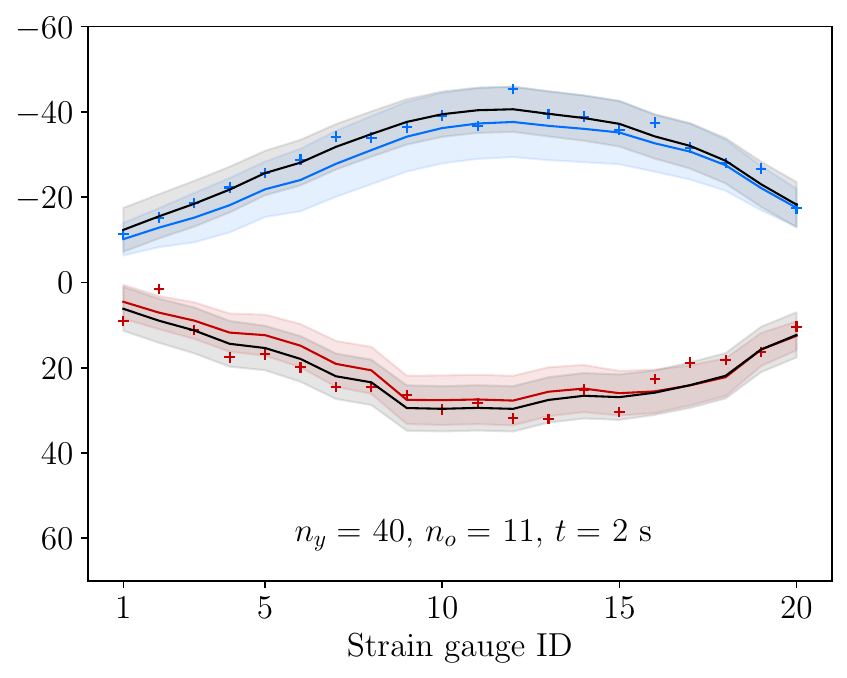}
	\includegraphics[width=0.31\textwidth]{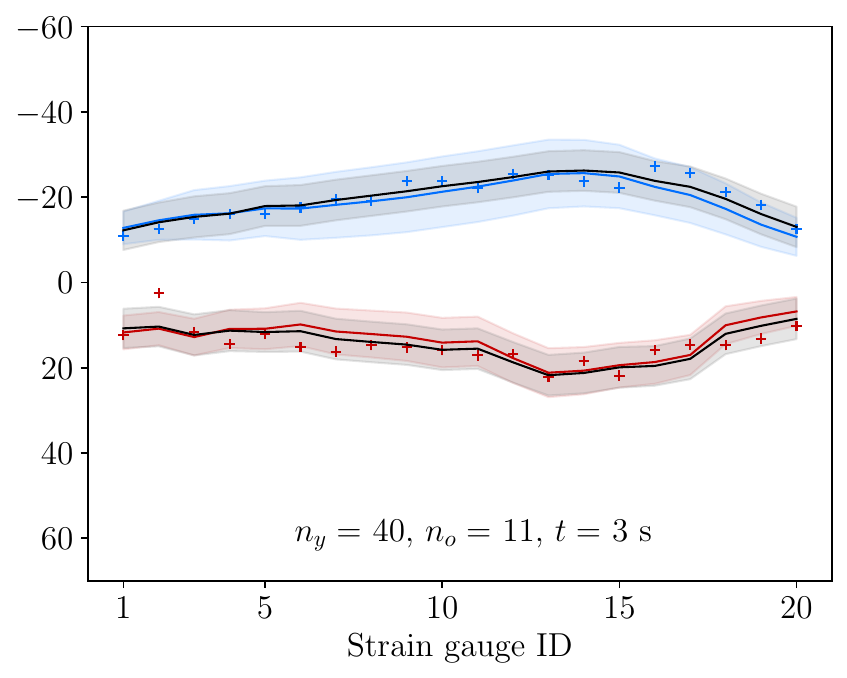}
	\caption{Inferred true strain density~$p(\ary z | \ary y)$  conditioned on strains (+) measured along the east main I-beam for~$n_o =11$ readings. The blue and red lines represent the mean~$\ary P \overline{ \ary u }$ of the prior along the top and bottom flanges, respectively, and the black lines the conditioned mean~$\overline {\ary z}_{|\ary y}$. The shaded areas denote the corresponding $95\%$ confidence regions. In each row the number of sensors~$n_y$ is fixed. In each column the observation time $t$ is fixed. The unit of the vertical axis is microstrain} \label{fig:pZgivYvarNy} 
\end{figure*}

\begin{figure*} []
	\setlength{\fboxsep}{0pt}%
	\setlength{\fboxrule}{0pt}%
	\centering
	\includegraphics[width=0.31\textwidth]{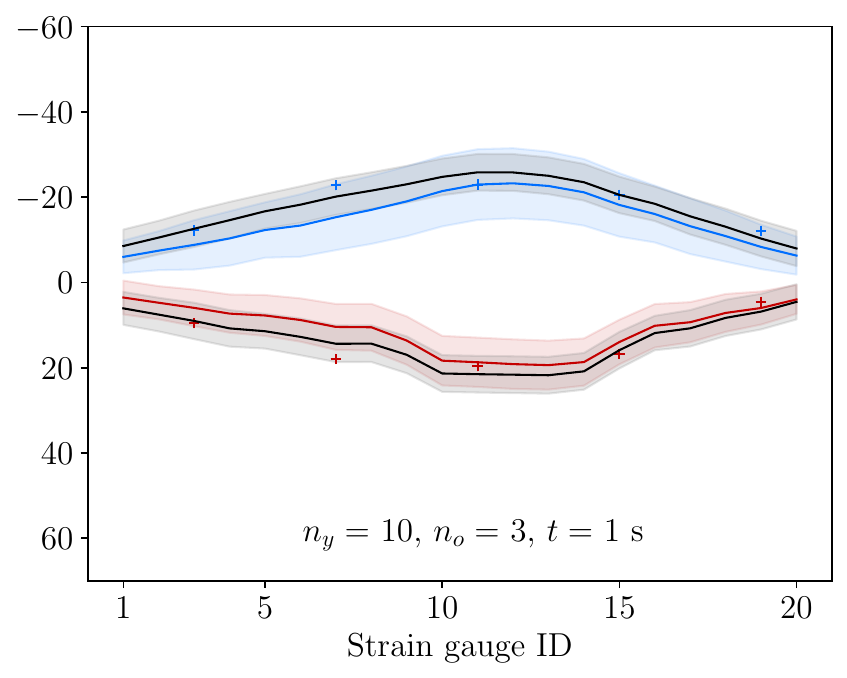} 
	\includegraphics[width=0.31\textwidth]{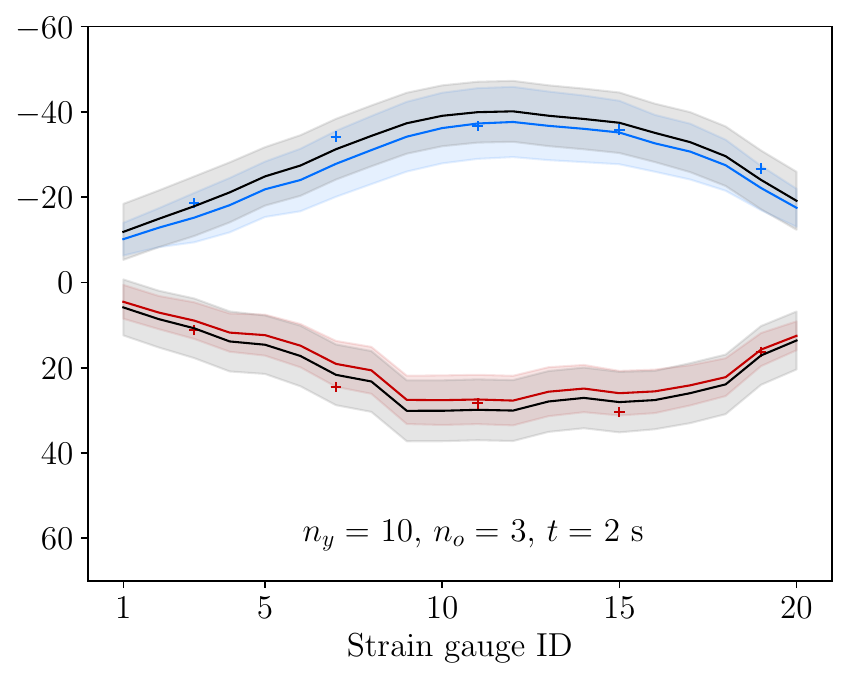}
	\includegraphics[width=0.31\textwidth]{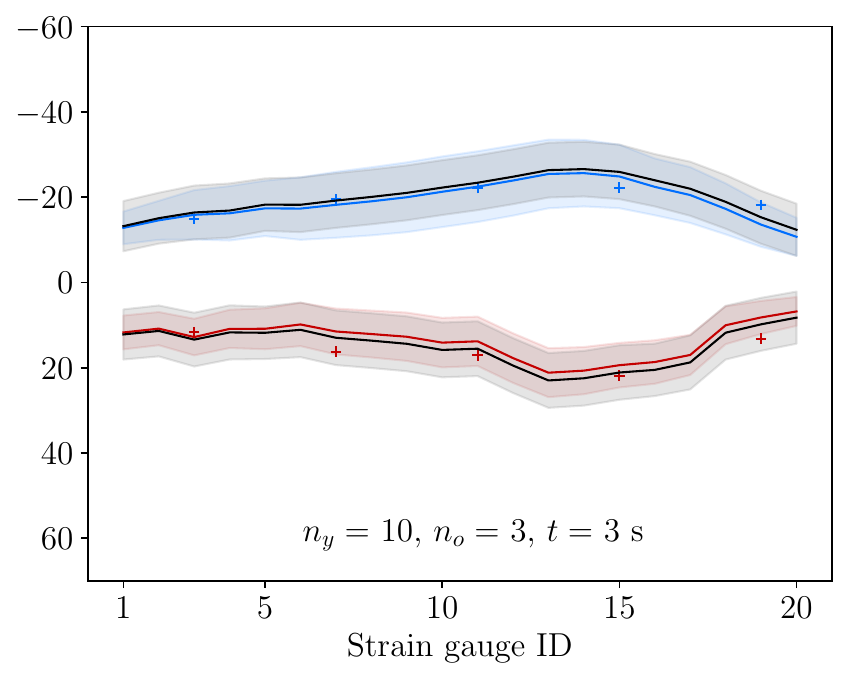}
	\\
	\includegraphics[width=0.31\textwidth]{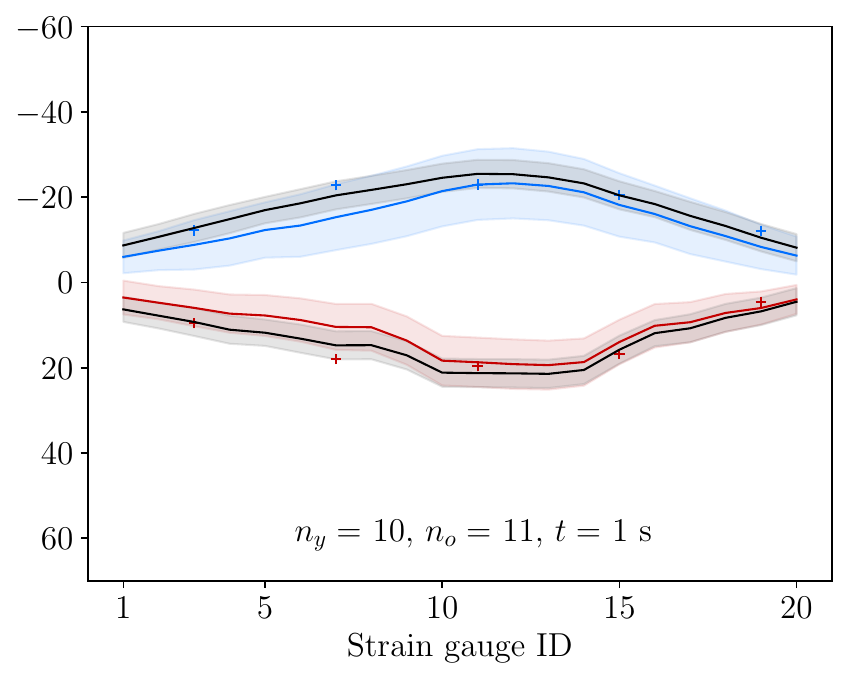} 
	\includegraphics[width=0.31\textwidth]{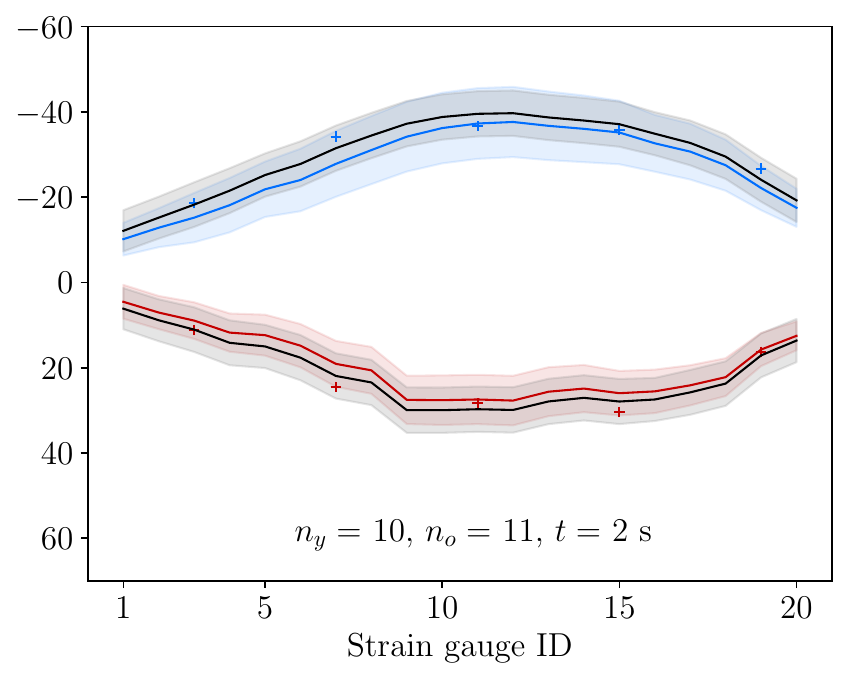}
	\includegraphics[width=0.31\textwidth]{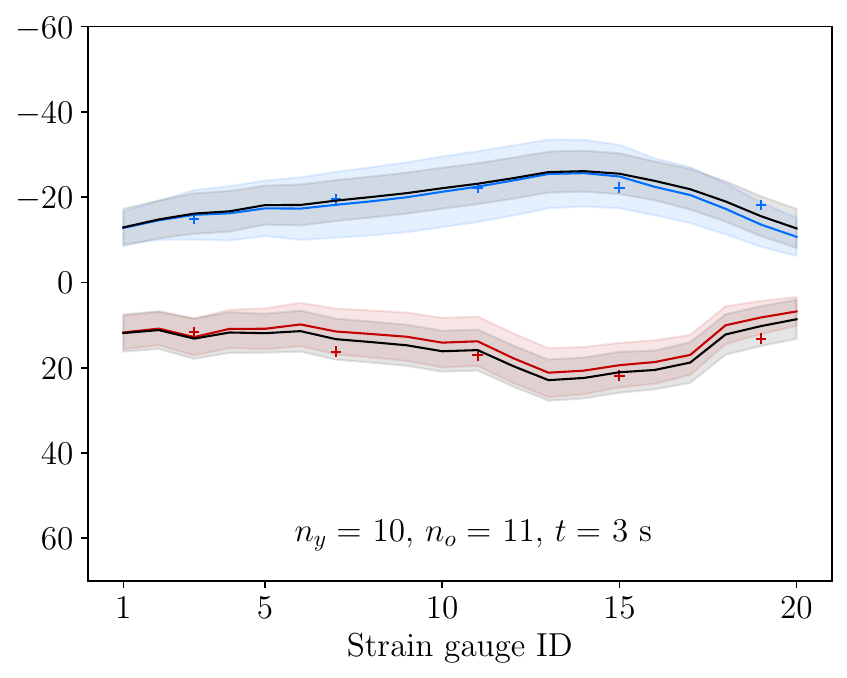}
	\\
	\includegraphics[width=0.31\textwidth]{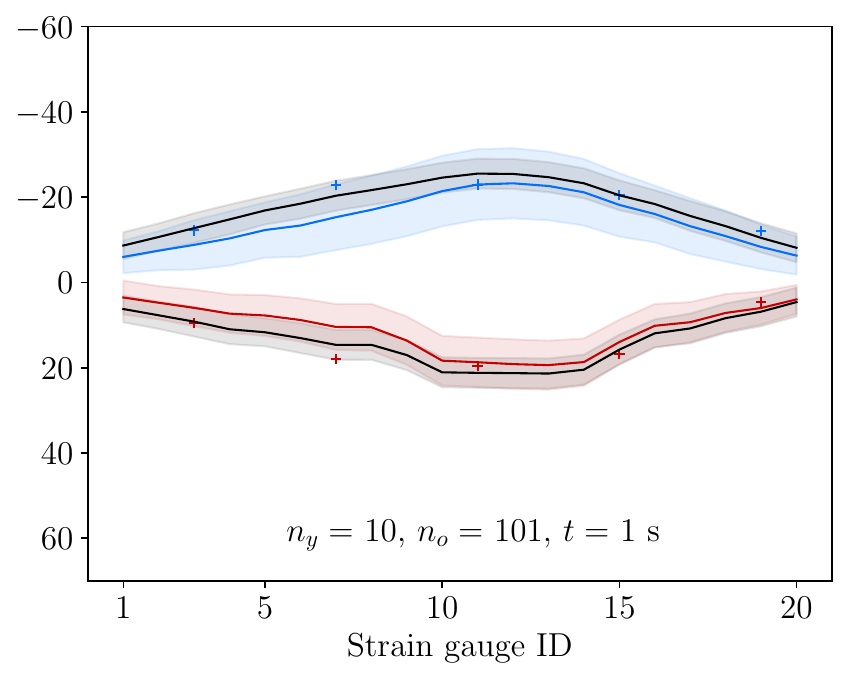} 
	\includegraphics[width=0.31\textwidth]{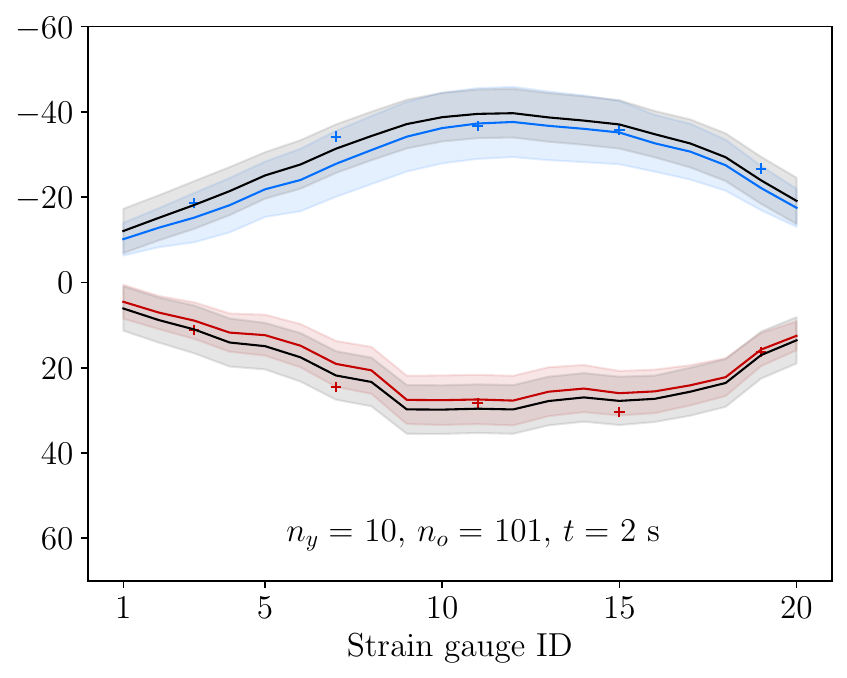}
	\includegraphics[width=0.31\textwidth]{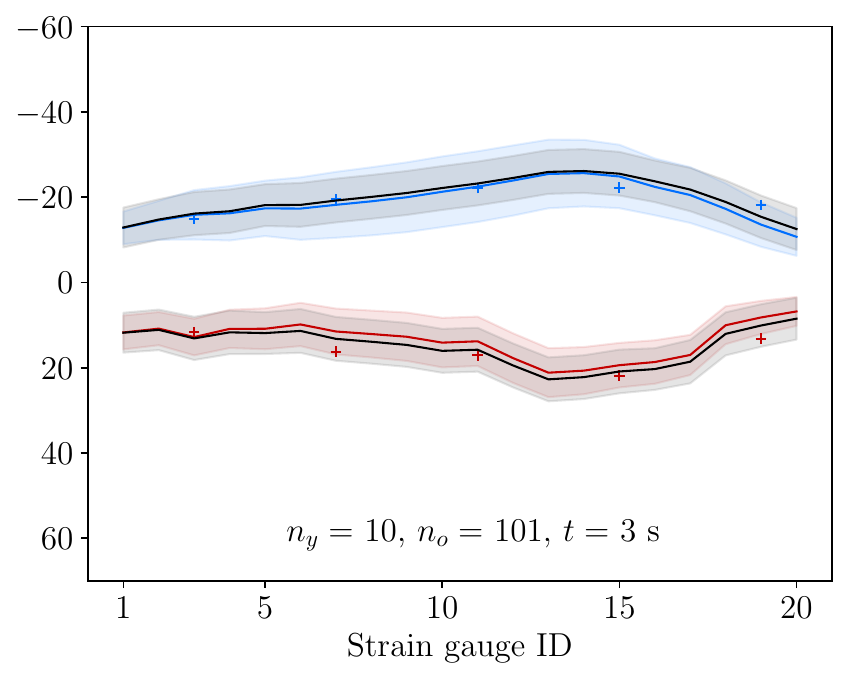}
	\\
	\includegraphics[width=0.31\textwidth]{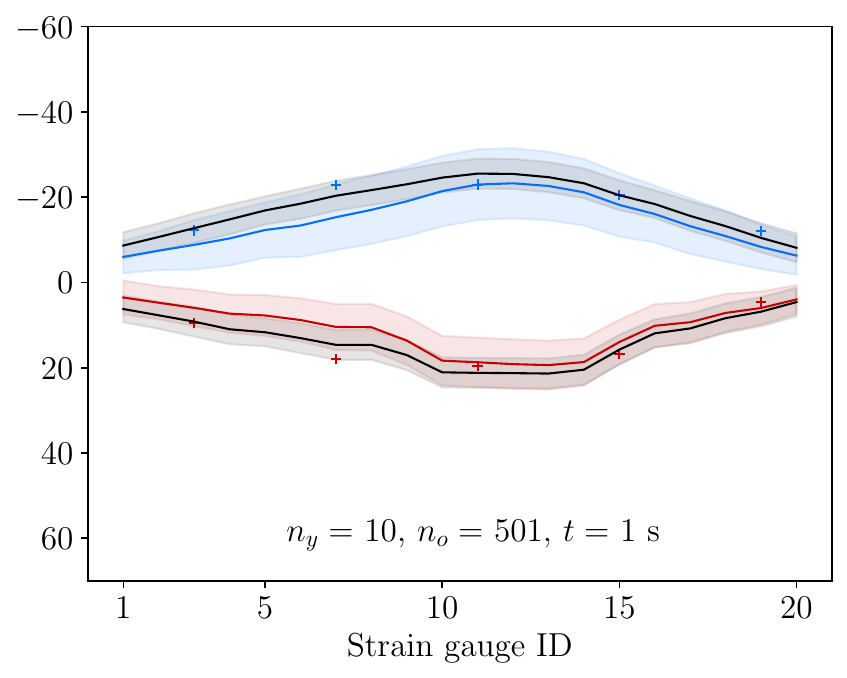} 
	\includegraphics[width=0.31\textwidth]{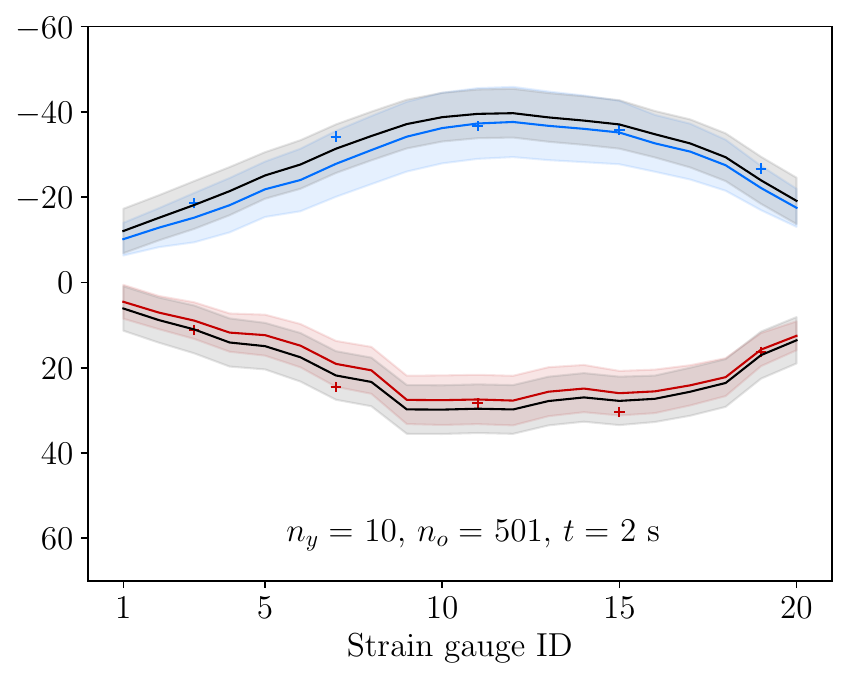}
	\includegraphics[width=0.31\textwidth]{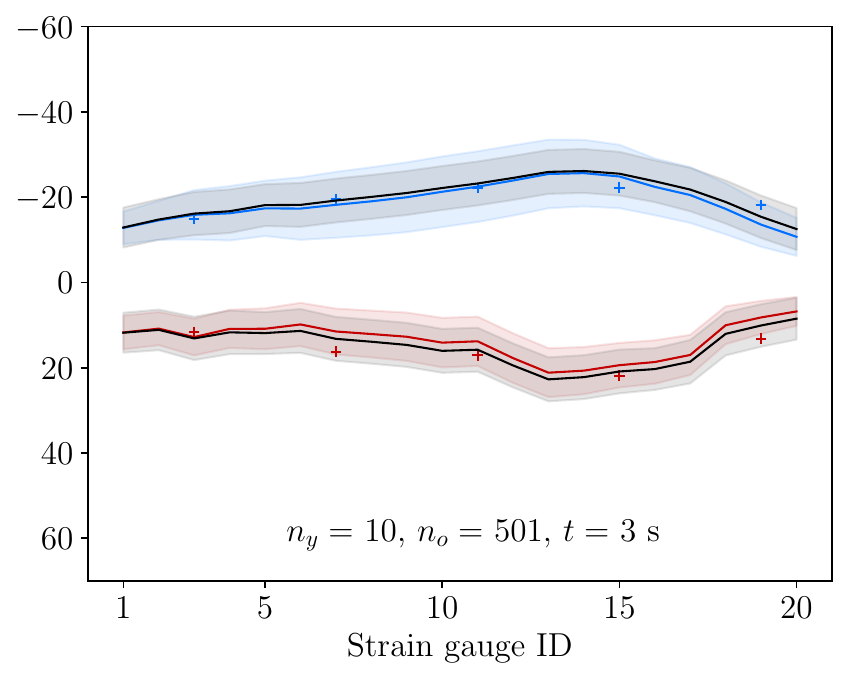}
	\caption{Inferred true strain density~$p(\ary z | \ary y)$  conditioned on strains (+) measured along the east main I-beam. The blue and red lines represent the mean~$\ary P \overline{\ary u }$ of the prior along the top and bottom flanges, respectively,  and the black lines the conditioned mean~$\overline {\ary z}_{| \ary y}$. The shaded areas denote the corresponding $95\%$ confidence regions. In each row the number of readings~$n_o$ and in each column the observation time $t$ is fixed. The unit of the vertical axis is microstrain} \label{fig:pZgivYvarNo} 
\end{figure*}

Next, the posterior true strain density~$p(\ary z | \ary y )$ is evaluated for different numbers of sensing points $n_y$ and a fixed number of readings \mbox{$n_o = 11$}, see Figure \ref{fig:pZgivYvarNy}. In contrast, Figure \ref{fig:pZgivYvarNo} shows the posterior true strains for different $n_o$ and a fixed $n_y = 10$. It can be observed that the means~$\overline{\ary z}_{| \ary y}$ obtained using fewer observation data show close agreement with those obtained using more data. Similarly, the 95\%  confidence intervals of the different posteriors $p(\ary z | \ary y)$ look visually very similar for all combinations of $n_o$ and $n_y$. 
These results confirm that it is possible to use data from fewer sensors $n_y$ and fewer readings  $n_o$ to obtain a sufficiently reliable estimate for the true strain response.

\subsubsection{Predictive strains at non-sensor locations \label{sec:westPrediction}}
%
Another advantage of statFEM is the possibility of generating strain predictions at locations where no sensing data is available. Limited sensing points may arise due to cost and labour considerations but it is also common in structural monitoring for sensor systems to malfunction and for entire sections of a sensor network to stop recording. The sensor strains along the west main I-beam are estimated by using only the readings from sensors installed along the east main I-beam. Although strain data along the west I-beam is available, they are not included in the observation matrix~$\ary Y$. This test is an extreme case of missing measurement data due either to sensor malfunction or temporary system error. 

Consider the vector of unobserved strains~$ \widehat{\ary y}$ at the locations with the coordinates~$\widehat{\ary X}$. The matrix~$\widehat{\ary X}$ contains the coordinates of the FBG sensors along the west I-beam, 
see Figure~\ref{fig:bridgeschematic}. The predictive density of sensors strains~$p( \widehat{\ary y}| \ary y)$ conditioned on the observations~$\ary y$ is given by
\begin{align} 
	p(\widehat{\ary y} | \ary y) = \int p(\widehat{\ary y} | \ary u) p(\ary u | \ary y) \D \ary u \, .
\end{align}
With the likelihood~\eqref{eq:likelihood} and the FE posterior~\eqref{eq:posteriorU} this becomes 
\begin{align} \label{eq:predictiveStrain}
	p(\widehat{\ary y} | \ary y) = \set N \left ( \rho \widehat{\ary P} \overline{\ary u}_{| \ary y}, \rho^2 \widehat{\ary P} \ary C_{\ary u | \ary y} \widehat{\ary P}^{\trans} + \widehat{\ary C}_{\ary d} + \widehat{\ary C}_{\ary e}  \right ) \, ,
\end{align}
where the matrices~$\widehat{\ary P}$, $\widehat{\ary C}_{\ary d}$ and $\widehat{\ary C}_{\ary e}$ are obtained by introducing the coordinates collected in~$\widehat{\ary X}$ in the respective operator and covariance kernels. Note that this expression is very similar to the true system density~$p(\ary z| \ary y)$  given in~\eqref{eq:posteriorZ} up to the additional covariance term due to the measurement errors. 

To evaluate the predictive strain density \mbox{$p(\widehat{\ary y} | \ary y)$}, the hyperparameters of the statistical model~$\ary w$ are learned first as before. To this end, included in the observation matrix $\ary Y$  are the strains of the $n_y = 40$ sensors along the east I-beam and for each sensor the $n_o = 101$ readings between the time \mbox{$1~\text{s} \leq t \leq 3~\text{s}$}. Subsequently, \mbox{$p(\widehat{\ary y} | \ary y)$} is computed by introducing the point estimate~\mbox{$\ary w^*$} obtained by MCMC sampling into~\eqref{eq:predictiveStrain}. The point estimate obtained in this example is \mbox{$\ary w^* = \begin{pmatrix}  0.87723 \,  &  3.996\,\text{microstrain} & 0.43 \, \text{m} \end{pmatrix}^\trans$}. The prediction $p(\widehat{\ary y} | \ary y)$ is computed at all the 40 sensor positions along the west I-beam. 

\begin{figure*} []
	\setlength{\fboxsep}{0pt}%
	\setlength{\fboxrule}{0pt}%
	\centering
	\includegraphics[width=0.31\textwidth]{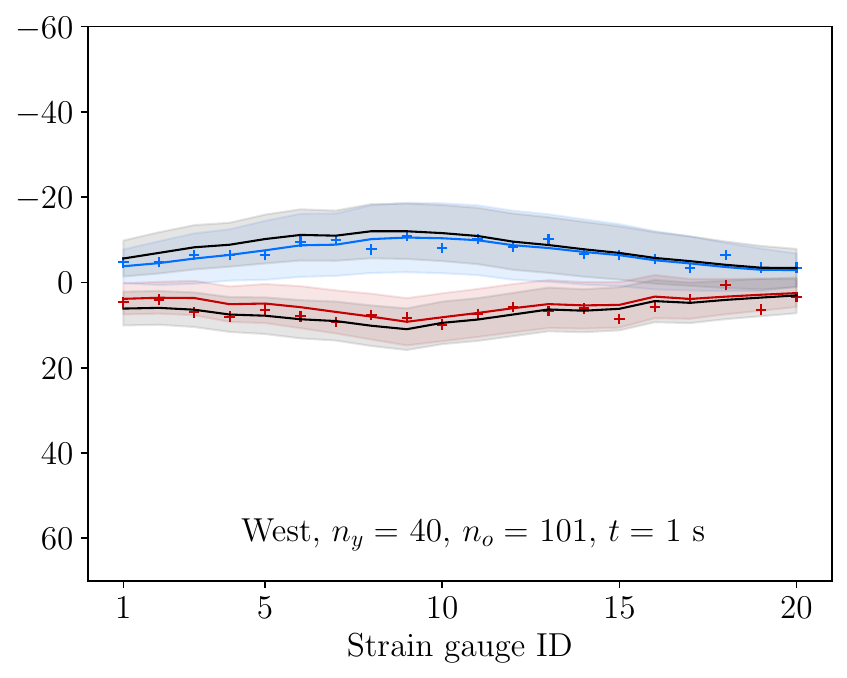}
	\includegraphics[width=0.31\textwidth]{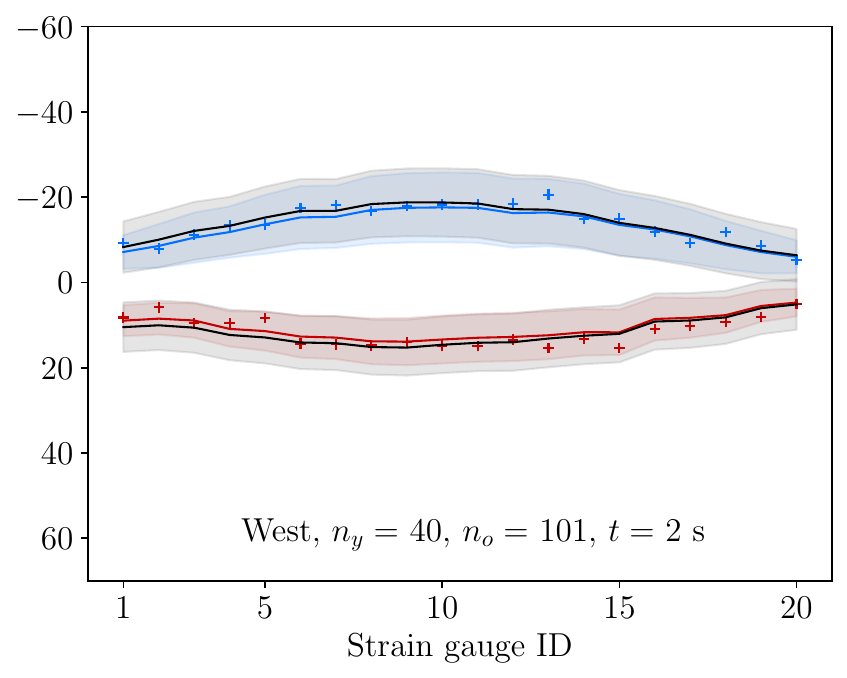}
	\includegraphics[width=0.31\textwidth]{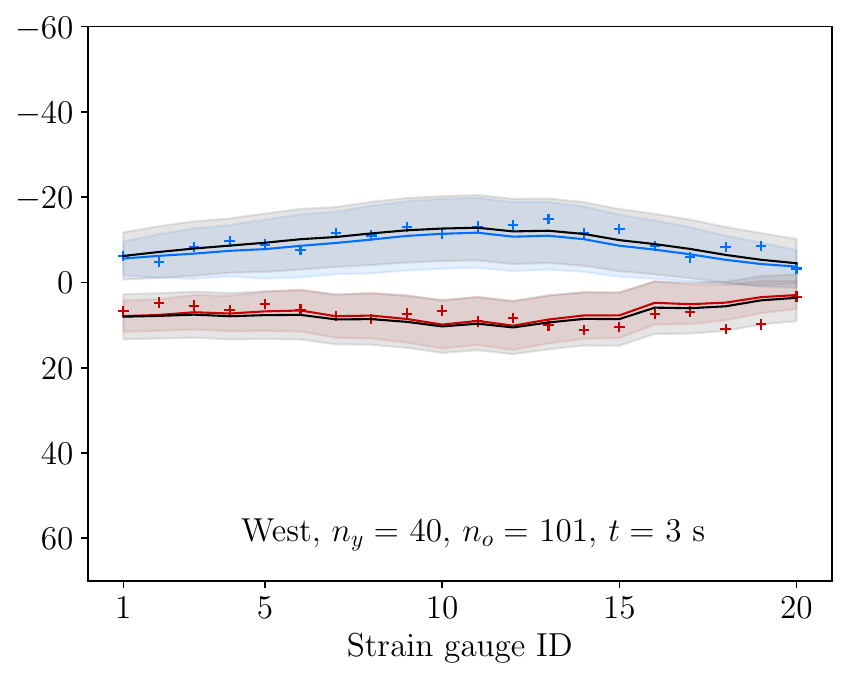}
	\caption{Predictive strain distribution  $p(\widehat{\ary y} | \ary y)$ along the west I-beam conditioned on strains (+) measured along the east main I-beam. 
	The blue and red lines represent the prior mean~$\ary P \overline{\ary u }$ along the top and bottom flanges, respectively, and the black lines the predictive mean~$\rho \ary P \overline{ \widehat{\ary u}}_{| \ary y}$. The shaded areas denote the corresponding $95\%$ confidence regions. The unit of the vertical axis is microstrain} \label{fig:pYgivYpred} 
\end{figure*}

As shown in Figure~\ref{fig:pYgivYpred}, with only measurement data from the east I-beam taken into account, the predictive strain distribution \mbox{$p(\widehat{\ary y} | \ary y)$} differs from the prior FE computed strains  $p(\ary P \ary u)$ and the measured strains. Its mean~$\overline{\widehat{\ary y}}_{\ary y}$ lies mostly in between the mean of the prior~$\ary P \overline{\ary u}$ and the measured strains. The 95\% confidence regions of the predictive strain distribution and the prior have almost the same width. Hence, even at locations where there is no measurement data available,  statFEM is able to improve the FE prior by utilising measurement data from other parts of the bridge. 

	%
\section{Conclusions \label{sec:conclusions}}
%

In this study, a statistical digital twin of an operational bridge superstructure has been developed to synthesise FBG strain measurements with FE model predictions for obtaining improved strain predictions. To this end, the application of statFEM as a new modelling paradigm in creating digital twins was proposed. 
In light of the presented results, several conclusions can be drawn.
\begin{enumerate}
	\item By considering strain data captured along the top and bottom flanges of the east and west main I-beams, statFEM analysis showed an improvement of the mean strain estimates after being conditioned on the recorded sensor data. The 95\% confidence bounds provided a quantifiable means of identifying anomalous sensor readings, which represents a significant step change in how structural monitoring data may be more reliably interpreted.

	\item To evaluate the effect of the number of sensing points on statFEM predictions, several sensor subsets ($40$, $20$ and $10$ sensors) for each main I-beam were considered. The resulting predictive strain distributions and associated 95\% confidence bounds between the three sensor subsets were negligible.  This demonstrated that statFEM may be used to optimise sensor network design leading to significant reductions in instrumentation costs.

	\item In addition to evaluating the effect of a reduced number of sensing points, several reduced sampling rates for each individual FBG sensor (originally capturing data at $250~\text{Hz}$) were considered. Once again, statFEM provided robust mean estimates of the predictive strains regardless of sampling rate.

	\item Using only strain data captured along the east main I-beam, the predictive strain distribution along the west main I-beam was estimated. It was shown that the predictive mean along the west I-beam incorporates the effects from both the FE prior and the measured data at the opposite east main I-beam.
	Therefore, in cases where missing or damaged parts of a sensor network exist, statFEM can still be used to generate meaningful interpretations of the data.
\end{enumerate}

In summary, this study has highlighted the suitability of statFEM for structural health monitoring and condition assessment of instrumented civil infrastructure assets. The presented results demonstrated that a relatively coarse FE model and a modest amount of data are sufficient to obtain reliable predictions. Combined with advanced computing technologies, this paves the way for future continuous, real-time monitoring and condition assessment of infrastructure assets. Furthermore, the choice of a suitable coarse FE model can be rationalised with the Bayesian model selection paradigm~\citep{beck2010bayesian, girolami2020statistical}. Indeed, Bayesian model selection  can provide a means to discover sudden changes in the state of the infrastructure by comparing the plausibility of different predefined `what-if' scenarios, including simulated damage and failure mechanisms.

\begin{Backmatter}
	
	\paragraph{Acknowledgments}
	We are grateful to the Cambridge Centre of Smart Infrastructure and Construction (CSIC) for providing access to the bridge monitoring data used in this study.
	
	\paragraph{Funding Statement}
	This work was supported by Wave 1 of The UKRI Strategic Priorities Fund under the EPSRC Grant EP/T001569/1, particularly the ``Digital twins for complex engineering systems" theme within that grant, and The Alan Turing Institute.
	
	\paragraph{Competing Interests}
	The authors declare no competing interests exist.
	
	\paragraph{Data Availability Statement}
	The data set used for this study is available via the University of Cambridge data repository.  The link to the repository can be found in~\cite{butler2018monitoring}.
	
	\paragraph{Author Contributions}
	Conseptualisation: E.F., L.B., M.G., and F.C.; Methodology: E.F., F.C.; Data curation: L.B.; Analysis, Validation, Visualisation: E.F., F.C.;  Writing original draft: E.F., L.B., F.C.; All authors approved the final submitted draft.
	
	\bibliographystyle{apalike-rv}
	\bibliography{smartBridge}

\end{Backmatter}	

\end{document}